\documentclass[11pt,oneside]{amsart}
\usepackage[T1]{fontenc}
\usepackage{geometry}
\usepackage{amssymb}

\makeatletter
 \theoremstyle{plain}
\newtheorem{thm}{Theorem}[section]
  \theoremstyle{plain}
  \newtheorem{cor}[thm]{Corollary}
  \theoremstyle{remark}
  \newtheorem{rem}[thm]{Remark}
  \theoremstyle{plain}
  \newtheorem{prop}[thm]{Proposition}
  \theoremstyle{plain}
  \newtheorem{lem}[thm]{Lemma}

\newcommand{\R}{\mathbb{R}}

\def\com#1{ \hbox{#1}}

\def\e{\hbox{\rm e}}

\smallskip
\def\<{{\langle }}
\def\>{{\rangle }}

\makeatother

\begin{document}

\title{Minimal tori with low nullity}

\author{David L. Johnson, Oscar Perdomo }

\date{\today}

\curraddr{Department of Mathematics\\
Lehigh University\\
Bethlehem, Pennsylvania 18015-3174}

\curraddr{Department of Mathematics\\
Central Connecticut State University\\
New Britain, CT 06050\\
}

\email{david.johnson@lehigh.edu, osperdom@univalle.edu.co}

\begin{abstract}
The \emph{nullity} of a minimal submanifold $M\subset S^{n}$ is the
dimension of the nullspace of the second variation of the area functional.
That space contains as a subspace the effect of the group of rigid
motions $SO(n+1)$ of the ambient space, modulo those motions which
preserve $M$, whose dimension is the \emph{Killing nullity} $kn(M)$
of $M$. In the case of 2-dimensional tori $M$ in $S^{3}$, there
is an additional naturally-defined 2-dimensional subspace; the dimension
of the sum of the action of the rigid motions and this space is the
\emph{natural nullity} $nnt(M)$. In this paper we will study minimal
tori in $S^{3}$ with natural nullity less than 8. We construct minimal
immersions of the plane $\R^{2}$ in $S^{3}$ that contain all possible
examples of tori with $nnt(M)<8$. We prove that the examples of Lawson
and Hsiang with $kn(M)=5$ also have $nnt(M)=5$, and we prove that
if the $nnt(M)\le6$ then the group of isometries of $M$ is not trivial. 
\end{abstract}


\maketitle

\section{Introduction}

Let $\tilde{\rho}:M\to S^{3}$ be a minimal immersion of an oriented
compact surface without boundary $M$ in the unit three dimensional
sphere $S^{3}\subset\R^{4}$. When there is no confusion, we will
identify $m\in M$ with $\tilde{\rho}(m)$ and the vectors in $T_{m}M$
with those in $d\tilde{\rho}_{m}(T_{m}M)\subset\R^{4}$. Let $N:M\to S^{3}$
be a Gauss map, i.e. $N(m)\perp T_{m}M$ and $\<N(m),{m}\>=0$. For
any $m\in M$, $a(m)$ will denote the nonnegative principal curvature
of $M$ at $m$ and $W_{1}(m)$ and $W_{2}(m)$ will denote two unit
tangent vectors such that $dN_{m}(W_{1}(m))=-a(m)W_{1}(m)$ and $dN_{m}(W_{2}(m))=a(m)W_{2}(m)$.
Notice that $a(m)$ is uniquely determined but the vectors $W_{1}(m)$
and $W_{2}(m)$ are not. When $M$ is a torus, it is known that for
every $m$, $a(m)$ is positive \cite{L}, therefore in this case
we can choose $W_{1}(m)$ and $W_{2}(m)$ so that they define smooth
vector field in $M$. In the following, if $M$ is a torus, $W_{1}$
and $W_{2}$ will denote such unit tangent vector fields and $a:M\to{\bf \R}$
will be the smooth function given by the positive principal curvature.
Since we are identifying vectors in $T_{m}M$ with those in $d\tilde{\rho}_{m}(T_{m}M)$
using the homomorphisms $d\tilde{\rho}_{m}$, the tangent vector fields
on $M$ are given by functions $X:M\to\R^{4}$ such that $\<X(m),\tilde{\rho}(m)\>=0$
and $\<X(m),N(m)\>=0$ for all $m\in M$. Given a fixed, skew-symmetric
$4\times4$ matrix $B$, define $f_{B}:M\to\R$ by $f_{B}=\<B\tilde{\rho}(m),N(m)\>.$
Since $M$ is minimal, $M$ is a critical point of the area functional.
The second variation of the area function at this critical point is
given by the stability operator\[
J:C^{\infty}(M)\to C^{\infty}(M)\quad\com{given\, by}\quad J(f):=-\Delta f-2a^{2}f-2f.\]
It is clear that $f_{B}$ satisfies the elliptic equation $J(f_{B})=0$
because, when we move the immersion $M$ by the group of isometries
$\e^{Bt}:S^{3}\to S^{3}$ we induce a family that leaves the area
and second fundamental form constant; $f_{B}$ is the function associated
with this family. The \emph{nullity} of a minimal surface is defined
as the dimension of the kernel of the operator $J$ and will be denoted
by $n(M)$. In \cite{LH}, Lawson and Hsiang classify all the minimal
surfaces that are invariant under under a 1-parametric group of isometries
in $S^{3}$. As they point out at the end of their paper, one way
to see this classification is the following: Define

\[
KS=\{f_{B}:B\in so(4)\};\quad\com{Killing\, nullity=dim}(KS)=kn(M).\]

We have that $kn(M)\le n(M)$ and in general the Killing nullity is
expected to be $6$ since the dimension of $so(4)$ is $6$. Then,
they classify all the examples of surfaces with $kn(M)<6$. These
examples turn out to be the totally geodesic spheres with Killing
nullity $3$, the Clifford tori with Killing nullity $4$ and a collection
of tori with Killing nullity $5$. It is known that when $M$ is a
torus, for any angle $\theta$, the function $h_{\theta}:M\to\R$
given by

\[
h_{\theta}=\cos(\theta)a^{-\frac{3}{2}}W_{1}(a)+\sin(\theta)a^{-\frac{3}{2}}W_{2}(a)\]

satisfies that $J(h_{\theta})=0$. Actually, this equation was the
starting point in the classification of all constant mean curvature
surface in $\R^{3}$ and all minimal tori in $S^{3}$ given by Pinkall
and Sterling in \cite{PS}. The study of the nullity of minimal tori
will be helpful in understanding their examples. When $M$ is a torus,
we can define the following space and the following integer

\[
NS=\{f_{B}+\lambda h_{\theta}:B\in so(4),\ \lambda,\theta\in\R,\},\quad\com{Natural\, nullity\, for\, tori:=dim}(NS):=nnt(M).\]

Clearly $NS$ is a subset of the kernel of $J$ and therefore $nnt(M)\le n(M)$.
The natural nullity is expected to be $8$ because of the $6$ parameters
of $so(4)$ and the two parameters $\lambda$ and $\theta$ in the
definition of the space $NS$. 

In this paper we study minimal tori with natural nullity for tori
less than 8. We construct, in Theorem (\ref{thm:solutions-of-system}),
minimal immersions of the plane $\R^{2}$ in $S^{3}$ that contain
\emph{all} possible examples of tori with $nnt(M)<8$, which is shown
in Theorem (\ref{thm:nnt<8}). We prove (Corollary (\ref{cor:nnt=00003D5-iff-HL}))
that the examples of Lawson and Hsiang with $kn(M)=5$ are the only
immersed tori satisfying $nnt(M)=5$, although the question of the
(total) nullity of these examples is not resolved. Finally, we show
in Theorem (\ref{thm:nnt<6-isometries}) that if $nnt(M)\le6$ then
the group of isometries is not trivial.

\section{Preliminaries}

In this section we will review some known results that will be used
later on. The first result has already been used in the introduction
in order to define the unit tangent smooth vector fields $W_{1}$
and $W_{2}$ in an immersed minimal torus of $S^{3}$.

\begin{thm}
\textbf{\emph{{[}Lawson, \cite{L}]}} If $M\subset S^{3}$ is a closed
minimal surface and $a:M\to\R$ denotes the nonnegative principal
curvature function, then $a$ is positive everywhere if and only if
$\chi(M)=0$.
\end{thm}
The next theorem also was used in the introduction in order to defined
define the natural nullity for tori. Even though it is a known result,
for completeness sake we will provide a proof at the end of this section.

\begin{thm}
\label{natural-solutions}If $M\subset S^{3}$ is a minimal immersed
torus, and $W_{1}:M\to S^{3}$ and $W_{2}:M\to S^{3}$ are unit vector
field that define the principal directions, then the functions\[
h_{0},h_{\frac{\pi}{2}}:M\to\R\ \hbox{given\, by}\ h_{0}(m)=a^{-\frac{3}{2}}W_{1}(a)\ \com{and}\ h_{\frac{\pi}{2}}=a^{-\frac{3}{2}}W_{2}(a)\]
satisfy\[
J(h_{0})=-\Delta h_{0}-2h_{0}-2a^{2}h_{0}=0=J(h_{\frac{\pi}{2}}).\]

\end{thm}
The following theorem will be used in section 4 to prove that one
subfamily of the family of examples of minimal immersion of the plane
in $S^{3}$ we have constructed corresponds to the Lawson-Hsiang examples.

\begin{thm}
\textbf{\emph{\label{thm:Ramanaham}{[}Ramanaham \cite{R}]}} Let
$\tilde{\rho}:M\to S^{3}$ be a minimal immersion from an oriented
compact surface. Suppose that $M$ admits a one parameter group of
isometries $\phi_{t}:M\to M$ with respect to the induced metric.
Then, there exists a one-parameter family of orientation preserving
isometries $\Phi_{t}$ of $S^{3}$ such that $\tilde{\rho}\circ\phi_{t}=\Phi_{t}\circ\tilde{\rho}$
for all $t\in\R$.
\end{thm}
The next theorem is a consequence of the uniformization theorem applied
to a minimal torus in $S^{3}$.

\begin{thm}
\label{uniformization}For every minimal immersion of a torus $\tilde{\rho}:M\to S^{3}$,
there exists a covering map $\tau:\R^{2}\to M$, a doubly periodic
conformal immersion $\rho:\R^{2}\to S^{3}$, a Gauss map $\nu:\R^{2}\to S^{3}$,
and a fixed angle $\alpha$, so that\[
\rho(u,v)=\tilde{\rho}(\tau(u,v)),\quad\nu(u,v)\perp\rho_{*}(T_{(u,v)}\R^{2}),\quad\nu(u,v)\perp\rho(u,v),\]
and\begin{eqnarray*}
\frac{\partial^{2}\rho}{\partial u^{2}} & = & -\frac{\partial r}{\partial u}\frac{\partial\rho}{\partial u}+\frac{\partial r}{\partial v}\frac{\partial\rho}{\partial v}+\cos(2\alpha)\nu-\e^{-2r}\rho\\
\frac{\partial^{2}\rho}{\partial v^{2}} & = & \frac{\partial r}{\partial u}\frac{\partial\rho}{\partial u}-\frac{\partial r}{\partial v}\frac{\partial\rho}{\partial v}-\cos(2\alpha)\nu-\e^{-2r}\rho\\
\frac{\partial^{2}\rho}{\partial u\partial v} & = & -\frac{\partial r}{\partial v}\frac{\partial\rho}{\partial u}-\frac{\partial r}{\partial u}\frac{\partial\rho}{\partial v}-\sin(2\alpha)\nu\\
\frac{\partial\nu}{\partial u} & = & \e^{2r}(-\cos(2\alpha)\frac{\partial\rho}{\partial u}+\sin(2\alpha)\frac{\partial\rho}{\partial v})\\
\frac{\partial\nu}{\partial v} & = & \e^{2r}(\sin(2\alpha)\frac{\partial\rho}{\partial u}+\cos(2\alpha)\frac{\partial\rho}{\partial v})\end{eqnarray*}
where $\e^{-2r}=\<\frac{\partial\rho}{\partial u},\frac{\partial\rho}{\partial u}\>=\<\frac{\partial\rho}{\partial v},\frac{\partial\rho}{\partial v}\>$.
Moreover, $\Delta r+2\sinh(2r)=0$.
\end{thm}
\begin{proof}
The idea of the proof is the following: the existence of the conformal
map $\rho$ and the covering $\tau$ follows from the uniformization
theorem, the existence of the constant $\alpha$ follows from the
fact that\[
f(z)=f(u+iv)=\<\frac{\partial^{2}\rho}{\partial u^{2}},\nu\>-i\ \<\frac{\partial^{2}\rho}{\partial u\partial v},\nu\>\]
is an analytic, doubly periodic function in the whole plane, and therefore
it is constant. Clearly this constant function $f$ is not identically
zero otherwise $M$ would be totally geodesic. By scaling the coordinates
$u$ and $v$ by a constant, we can make $f(u+iv)=\cos(2\alpha)+i\sin(2\alpha)$
for some constant angle $\alpha$. 

To complete the proof, the equations for the second derivatives of
$\rho$ are just the standard computation of the Christoffel symbols
and the elliptic equation of $r$ follows from computing the Gauss
curvature using the Christoffel symbols and making it equal to $1-\e^{4r}$,
i.e, this elliptic equation follows from the Gauss equation. 
\end{proof}
\begin{cor}
\label{cor:prin-dirs}Using the same notation as in the previous theorem,
the principal directions of the minimal immersion are given by\[
V_{1}=\e^{r}(\cos(\alpha)\frac{\partial\rho}{\partial u}-\sin(\alpha)\frac{\partial\rho}{\partial v})\quad\com{and}\quad V_{2}=\e^{r}(\sin(\alpha)\frac{\partial\rho}{\partial u}+\cos(\alpha)\frac{\partial\rho}{\partial v}).\]
More precisely, \[
d\nu(\{d\rho_{(u,v)}\}^{-1}(W_{1}\circ\tau))=-\e^{2r}V_{1}\quad\com{and}\quad d\nu(\{d\rho_{(u,v)}\}^{-1}(W_{2}\circ\tau))=\e^{2r}V_{2}.\]
Moreover, it follows from the last expression that the principal curvatures
are $\pm a$ where the function $a:M\to\R$ satisfies $a(\tau(u,v))=\e^{2r(u,v)}$.
\end{cor}
\begin{rem}
\label{rem:rotatecoord}A direct computation shows that if make a
rotation of the variable $u$ and $v$, i.e. if we consider the variables
$\bar{u}$ and $\bar{v}$ where\[
u=\cos(\beta)\bar{u}+\sin(\beta)\bar{v}\quad\com{and}\quad v=-\sin(\beta)\bar{u}+\cos(\beta)\bar{v},\]
then, the angle $\alpha$ in the theorem above changes from $\alpha$
to $\alpha-\beta$. 
\end{rem}
\begin{cor}
\label{cor:a}If $M\subset S^{3}$ is a minimal immersed torus, $N:M\to S^{3}$
its Gauss map, and $W_{1}:M\to S^{3}$ and $W_{2}:M\to S^{3}$ are
unit vector field that define the principal directions with $dN_{m}(W_{1})=-aW_{1}$
and $dN_{m}(W_{2})=aW_{2}$, where $a:M\to\R$ is the positive principal
curvature function, then\begin{eqnarray*}
\bar{\nabla}_{W_{1}}W_{1} & = & \frac{W_{2}(a)}{2a}W_{2}+aN-m\\
\bar{\nabla}_{W_{1}}W_{2} & = & -\frac{W_{2}(a)}{2a}W_{1}=\nabla_{W_{1}}W_{2}\\
\bar{\nabla}_{W_{2}}W_{1} & = & -\frac{W_{1}(a)}{2a}W_{2}=\nabla_{W_{2}}W_{1}\\
\bar{\nabla}_{W_{2}}W_{2} & = & \frac{W_{1}(a)}{2a}W_{1}-aN-m\end{eqnarray*}
where $\bar{\nabla}$ is the Levi-Civita Connection in $\R^{4}$ and
$\nabla$ is the connection in the surface $M$.
\end{cor}
\begin{proof}
This follows from Theorem (\ref{uniformization}) and Corollary (\ref{cor:prin-dirs})
by noticing that

\[
W_{1}(\tau(u,v))=V_{1}(u,v),\quad W_{2}(\tau(u,v))=V_{2}(u,v)\quad\hbox{and}\quad N(\tau(u,v))=\nu(u,v).\]

\end{proof}
\begin{rem}
Theorem (\ref{natural-solutions}) can be proved using Theorem (\ref{uniformization})
in the following way: Take maps $\rho,V_{1},V_{2},\nu:\R^{2}\to S^{3}$,
$\tau:\R^{2}\to M$ and $r:\R^{2}\to\R$ such that they satisfy the
condition of Theorem (\ref{uniformization}) with $\alpha=0$, i.e.,
with $V_{1}(u,v)=W_{1}(\tau(u,v))=\e^{r(u,v)}\frac{\partial\rho}{\partial u}(u,v)$
and $V_{2}(u,v)=W_{2}(\tau(u,v))=\e^{r(u,v)}\frac{\partial\rho}{\partial v}(u,v)$.
Since $\Delta_{\R^{2}}r+2\sinh{(2r)}=0$ we obtain that
\end{rem}
\[
\Delta_{\R^{2}}\frac{\partial r}{\partial u}+4\cosh{(2r)}\frac{\partial r}{\partial u}=0.\]

Since $\frac{\partial\rho}{\partial u}(u,v)=\e^{-r}V_{1}(u,v)=\e^{-r}W_{1}(\tau(u,v))$
and $a(\tau(u,v))=\e^{2r(u,v)}$, we have

\[
\frac{\partial r}{\partial u}=a^{-\frac{1}{2}}W_{1}(\frac{1}{2}\ln(a))=\frac{1}{2}a^{-\frac{3}{2}}W_{1}(a).\]

Denote by $\Delta_{M}$ the Laplacian in the surface. Since the metric
induced by $\rho$ in $\R^{2}$ is given by $ds^{2}=\e^{-2r}(du^{2}+dv^{2})$,
we obtain that,

\[
\Delta_{M}(\frac{1}{2}a^{-\frac{3}{2}}W_{1}(a))=a\Delta_{\R^{2}}(\frac{\partial r}{\partial u})=-a(2(a+a^{-1})(\frac{1}{2}a^{-\frac{3}{2}}W_{1}(a)))\]

Therefore the function $h_{0}=a^{-\frac{3}{2}}W_{1}(a)$ satisfies
$J(h_{0})=0$. We prove that $J(h_{\frac{\pi}{2}})=0$ similarly.

\section{Integrable systems and solutions of the sinh-Gordon Equation}

In this section we will study an integrable system that produces solutions
of the sinh-Gordon equations, the construction made here is similar
to that of \cite{Ab} and \cite{PS}. The integrable system lives
in

\[
\R^{{18}}=\{(p,V_{1},V_{2},\nu,r,s):p,V_{1},V_{2},\nu\in\R^{4},\quad\hbox{and}\quad r,s\in\R\}\]

and is given by the vector fields $Z,W:\R^{{18}}\to\R^{{18}}$ given
by\begin{eqnarray*}
Z & = & \big{(}\e^{-r}(\cos(\theta)V_{1}+\sin(\theta)V_{2}),sV_{2}+\cos(\theta)(\e^{r}\nu-\e^{-r}p),\\
 &  & -sV_{1}-\sin(\theta)(\e^{r}\nu+\e^{-r}p),\e^{r}(-\cos(\theta)V_{1}+\sin(\theta)V_{2}),\\
 &  & \langle Bp,\nu\rangle,\cos(\theta)\<BV_{2},\e^{-r}\nu-\e^{r}p\>-\sin(\theta)\<BV_{1},\e^{-r}\nu+\e^{r}p\>\big{)}\\
 & = & (Z_{1},Z_{2},Z_{3},Z_{4},Z_{5},Z_{6})\\
\\W & = & \big{(}\e^{-r}(\cos(\theta)V_{2}-\sin(\theta)V_{1}),-\<Bp,\nu\>V_{2}+\sin(\theta)(-\e^{r}\nu+\e^{-r}p),\\
 &  & \<Bp,\nu\>V_{1}-\cos(\theta)(\e^{r}\nu+\e^{-r}p),\e^{r}(\sin(\theta)V_{1}+\cos(\theta)V_{2}),s,\\
 &  & \e^{-2r}-\e^{2r}-\sin(\theta)\<BV_{2},\e^{-r}\nu-\e^{r}p\>-\cos(\theta)\<BV_{1},\e^{-r}\nu+\e^{r}p\>\big{)}\\
 & = & (W_{1},W_{2},W_{3},W_{4},W_{5},W_{6}),\end{eqnarray*}

where $B$ is a skew-symmetric $4\times4$ matrix, i.e. $B\in so(4)$,
and $\theta$ is any real number. In the notation above, $W_{i}$
and $Z_{i}$ have values in $\R^{4}$ for $i=1,2,3,4$, and in $\R$
for $i=5,6$. We can write the system induced by these vector fields
as\begin{eqnarray*}
\frac{\partial p}{\partial u} & = & \e^{-r}(\cos(\theta)V_{1}+\sin(\theta)V_{2})\\
\frac{\partial V_{1}}{\partial u} & = & sV_{2}+\cos(\theta)(\e^{r}\nu-\e^{-r}p)\\
\frac{\partial V_{2}}{\partial u} & = & -sV_{1}-\sin(\theta)(\e^{r}\nu+\e^{-r}p)\\
\frac{\partial\nu}{\partial u} & = & \e^{r}(-\cos(\theta)V_{1}+\sin(\theta)V_{2})\\
\frac{\partial r}{\partial u} & = & \langle Bp,\nu\rangle\\
\frac{\partial s}{\partial u} & = & \cos(\theta)\<BV_{2},\e^{-r}\nu-\e^{r}p\>-\sin(\theta)\<BV_{1},\e^{-r}\nu+\e^{r}p\>\end{eqnarray*}
and\begin{eqnarray}
\frac{\partial p}{\partial v} & = & \e^{-r}(\cos(\theta)V_{2}-\sin(\theta)V_{1})\nonumber \\
\frac{\partial V_{1}}{\partial v} & = & -\<Bp,\nu\>V_{2}+\sin(\theta)(-\e^{r}\nu+\e^{-r}p)\nonumber \\
\frac{\partial V_{2}}{\partial v} & = & \<Bp,\nu\>V_{1}-\cos(\theta)(\e^{r}\nu+\e^{-r}p)\label{eq:1}\\
\frac{\partial\nu}{\partial v} & = & \e^{r}(\sin(\theta)V_{1}+\cos(\theta)V_{2})\nonumber \\
\frac{\partial r}{\partial v} & = & s\nonumber \\
\frac{\partial s}{\partial v} & = & \e^{-2r}-e^{2r}-\sin(\theta)\<BV_{2},\e^{-r}\nu-\e^{r}p\>-\cos(\theta)\<BV_{1},\e^{-r}\nu+\e^{r}p\>\nonumber \end{eqnarray}

We will refer to the previous system as the \emph{integrable system}
(\ref{eq:1}).

The following theorem provides a family a solutions of the sinh-Gordon
equation.

\begin{thm}
The vector fields $Z$ and $W$ commute, and if\[
\Theta_{Z}:(-\epsilon,\epsilon)\times\R^{{18}}\to\R^{{18}}\quad\hbox{and}\quad\Theta_{W}:(-\epsilon,\epsilon)\times\R^{{18}}\to\R^{{18}}\]
are the flows of the vector fields $Z$ and $W$ respectively, and
for any ${\bf x_{0}}\in\R^{{18}}$, we define the map $\phi:\R^{2}\to\R^{{18}}$
by\[
\phi(u,v)=\Theta_{Z}(u,\Theta_{W}(v,{\bf x^{0}}))=(\phi_{1}(u,v),\dots,\phi_{18}(u,v)),\]
then, the function $r(u,v)=\phi_{17}(u,v)$ solves the equation\[
\Delta r+2\sinh(2r)=0.\]

\end{thm}
\begin{proof}
We will denote by $DW$ and $DZ$ the $18\times18$ matrices of the
first derivatives of $W$ and $Z$ respectively. We will also denote
by $V(f)$ the directional derivative of the function $f$ in the
direction of $V$, the function $f$ may be a vector value function.
For example\[
Z_{5}(f)=\<Bp,\nu\>\frac{\partial}{\partial r}(f)\ \com{and}\ W_{1}(f)=\sum_{i=1}^{4}\e^{-r}(\cos(\theta)V_{2}^{i}-\sin(\theta)V_{2}^{i})\frac{\partial}{\partial p^{i}}(f)\]
Notice that\[
Z_{1}(p)=Z_{1},\ Z_{1}(V_{i})=0,\ Z_{1}(\nu)=0,\ Z_{1}(\e^{\pm r})=0,\ Z_{1}(s)=0\]
We get similar equations for $Z_{i}(p)$, $Z_{i}(V_{j})$, $Z_{i}(\nu)$,
$Z_{i}(\e^{\pm r})$, $Z_{i}(s)$, $W_{i}(p)$, $W_{i}(V_{j})$, $W_{i}(\nu)$,
$W_{i}(\e^{\pm r})$, and $W_{i}(s)$. 

Now,\begin{eqnarray*}
[Z,W] & = & (DW)Z-(DZ)W=ZW-WZ\\
 & = & (Z(W_{1}),Z(W_{2}),Z(W_{3}),Z(W_{4}),Z(W_{5}),Z(W_{6}))-\\
 &  & (W(Z_{1}),W(Z_{2}),W(Z_{3}),W(Z_{4}),W(Z_{5}),W(Z_{6}))\end{eqnarray*}
The first four components of the vector above are given by\begin{eqnarray*}
Z(W_{1})-W(Z_{1}) & = & \e^{-r}\cos(\theta)Z_{3}-\e^{-r}\sin(\theta)Z_{2}-\e^{-r}(\cos(\theta)V_{2}-\sin(\theta)V_{1})Z_{5}\\
 &  & -(\e^{-r}\cos(\theta)W_{2}+\e^{-r}\sin(\theta)W_{3}-\e^{-r}(\cos(\theta)V_{1}+\sin(\theta)V_{2})W_{5})\\
 & = & \e^{-r}\{\cos(\theta)(-sV_{1}-\sin(\theta)(\e^{r}\nu+\e^{-r}p))-\sin(\theta)(sV_{2}\\
 &  & +\cos(\theta)(-\e^{-r}p+\e^{r}\nu))\}-\<Bp,\nu\>\e^{-r}(\cos(\theta)V_{2}-\sin(\theta)V_{1})\\
 &  & -\e^{-r}(\cos(\theta)(-\<Bp,\nu\>V_{2}+\sin(\theta)(\e^{-r}p-\e^{r}\nu))+\sin(\theta)(\<Bp,\nu\>V_{1}\\
 &  & -\cos(\theta)(\e^{-r}p+\e^{r}\nu)))+\e^{-r}(\cos(\theta)V_{1}+\sin(\theta)V_{2})s\\
 & = & 0.\end{eqnarray*}
Direct computations, some of them longer, some others shorter than
the above, show that the other components of $[Z,W]$ are also zero.
We now show that $r(u,v)=\phi_{17}(u,v)$ is a solution of the sinh-Gordon
equation. We have that\begin{eqnarray*}
\Delta r & = & \frac{\partial^{2}r}{\partial u^{2}}+\frac{\partial^{2}r}{\partial v^{2}}=\frac{\partial\ \<Bp,\nu\>}{\partial u}+\frac{\partial s}{\partial v}\\
 & = & \<B(\e^{-r}(\cos(\theta)V_{1}+\sin(\theta)V_{2})),\nu\>+\<Bp,\e^{r}(-\cos(\theta)V_{1}+\sin(\theta)V_{2})\>\\
 &  & -2\sinh(2r)-\sin(\theta)\<BV_{2},\e^{-r}\nu-\e^{r}p\>-\cos(\theta)\<BV_{1},\e^{-r}\nu+\e^{r}p\>\\
 & = & -2\sinh(2r).\end{eqnarray*}
Notice that in the last step we have used the fact that $\<Bp,V_{2}\>=-\<BV_{2},p\>$,
i.e. we have used the fact that $B^{T}=-B$. 
\end{proof}
The previous theorem shows that for any choice of $B\in so(4)$, $\theta\in\R$
and ${\bf x}_{0}\in\R^{{18}}$ we have a solution of the sinh-Gordon
equation. One may think that a similar integrable system in $\R^{{4n+2}}=\R^{n}\times\R^{n}\times\R^{n}\times\R^{n}\times\R^{2}$
can be defined so that it produces a bigger space of solutions for
the sinh-Gordon equation. Indeed the system (\ref{eq:1}) can be generalized
to an integrable system in $\R^{{4n+2}}$ for any given $B\in so(n)$
and $\theta\in\R$, but the solutions of the sinh-Gordon will reduce
to solutions in the case $n=1$, as the following proposition explains.

\begin{prop}
If we think about the integrable system (\ref{eq:1}) as being defined
in $\R^{{4n+2}}$ by taking the vectors $p,V_{1},V_{2}$ and $\nu$
in $\R^{n}$ instead of vectors in $\R^{4}$, and by taking a skew
symmetry matrix $B\in so(n)$ instead of a matrix in $so(4)$, then\renewcommand{\labelenumi}{(\alph{enumi})}
\begin{enumerate}
\item The new system is integrable.
\item If $\phi:(-\epsilon,\epsilon)\times(-\epsilon,\epsilon)\to\R^{{4n+2}}$
is the solution of the new system with initial condition $x_{0}=(p^{0},V_{1}^{0},V_{2}^{0},\nu^{0},(r^{0},s^{0}))\in\R^{{4n+2}}$,
then
\end{enumerate}
\begin{eqnarray*}
\phi(u,v) & = & \ (\tilde{\phi}_{1}\tilde{e}_{1}+\tilde{\phi}_{2}\tilde{e}_{2}+\tilde{\phi}_{3}\tilde{e}_{3}+\tilde{\phi}_{4}\tilde{e}_{4},\ \tilde{\phi}_{5}\tilde{e}_{1}+\tilde{\phi}_{6}\tilde{e}_{2}+\tilde{\phi}_{7}\tilde{e}_{3}+\tilde{\phi}_{8}\tilde{e}_{4},\ \tilde{\phi}_{9}\tilde{e}_{1}+\tilde{\phi}_{10}\tilde{e}_{2}+\\
 &  & \tilde{\phi}_{11}\tilde{e}_{3}+\tilde{\phi}_{12}\tilde{e}_{4},\ \tilde{\phi}_{13}\tilde{e}_{1}+\tilde{\phi}_{14}\tilde{e}_{2}+\tilde{\phi}_{15}\tilde{e}_{3}+\tilde{\phi}_{16}\tilde{e}_{4},(\tilde{\phi}_{17},\tilde{\phi}_{18}))\end{eqnarray*}
where $\{\tilde{e}_{1},\tilde{e}_{2},\tilde{e}_{3},\tilde{e}_{4}\}$
is an orthonormal basis of any 4 dimensional subspace $S$ of $\R^{n}$
that contains the vectors $p^{0},V_{1}^{0},V_{2}^{0}$ and $\nu^{0}$,
and $\tilde{\phi}:(-\epsilon,\epsilon)\times(-\epsilon,\epsilon)\to\R^{{18}}$
is the solution of the integrable system (\ref{eq:1}) with initial
condition $(\tilde{p}^{0},\tilde{V}_{1}^{0},\tilde{V}_{2}^{0},\tilde{\nu}^{0},(r^{0},s^{0}))$,
where $\tilde{p}^{0},\tilde{V}_{1}^{0},\tilde{V}_{2}^{0},\tilde{\nu}^{0}$
are the coordinates of the vectors $p^{0},V_{1}^{0},V_{2}^{0},\nu^{0}$
in $S$ with respect to the basis $\{\tilde{e}_{1},\tilde{e}_{2},\tilde{e}_{3},\tilde{e}_{4}\}$,
and matrix $\tilde{B}\in so(4)$ given by $\tilde{b}_{ij}=\<\tilde{B}\tilde{e}_{i},\tilde{e}_{j}\>$.
\end{prop}
\begin{proof}
The proof of part (a) follows the same lines as the first part in
the proof of Theorem 3.1. Part (b) follows form the fact that, if
for any $y\in\R^{4}$ we define $\hat{y}=y_{1}\tilde{e}_{1}+y_{2}\tilde{e}_{2}+y_{3}\tilde{e}_{3}+y_{4}\tilde{e}_{4}$,
then $\<B\hat{y},\hat{z}\>=\<\tilde{B}y,z\>$ and $\widehat{F^{\prime}(u)}=\hat{F}^{\prime}(u)$
for any differentiable map $F:\R\to\R^{4}$. 
\end{proof}
In order to study the system (\ref{eq:1}) we will define this second
integrable system,

\begin{thm}
\label{thm:system-2}Let $p,\nu,V_{1},V_{2}:(-\epsilon,\epsilon)\times(-\epsilon,\epsilon)\to\R^{4}$
and $r,s:(-\epsilon,\epsilon)\times(-\epsilon,\epsilon)\to\R$ be
a solution of the integrable system (\ref{eq:1}). If $\tilde{B}$
is another skew symmetric matrix, and if we define the functions\begin{eqnarray*}
\xi_{1} & = & \<Bp,\nu\>,\ \xi_{2}=\<BV_{1},p\>,\ \xi_{3}=\<BV_{2},p\>,\ \xi_{4}=\<BV_{1},V_{2}\>,\ \xi_{5}=\<BV_{1},\nu\>,\ \xi_{6}=\<BV_{2},\nu\>\\
\tilde{\xi}_{1} & = & \<\tilde{B}p,\nu\>,\ \tilde{\xi}_{2}=\<\tilde{B}V_{1},p\>,\ \tilde{\xi}_{3}=\<\tilde{B}V_{2},p\>,\ \tilde{\xi}_{4}=\<\tilde{B}V_{1},V_{2}\>,\ \tilde{\xi}_{5}=\<\tilde{B}V_{1},\nu\>,\ \tilde{\xi}_{6}=\<\tilde{B}V_{2},\nu\>,\end{eqnarray*}
then\begin{eqnarray*}
\frac{\partial\xi_{1}}{\partial u} & = & \ \e^{r}(\cos(\theta)\xi_{2}-\sin(\theta)\xi_{3})+\e^{-r}(\cos(\theta)\xi_{5}+\sin(\theta)\xi_{6})\\
\frac{\partial\xi_{2}}{\partial u} & = & \ s\xi_{3}-\e^{r}\cos(\theta)\xi_{1}+\e^{-r}\sin(\theta)\xi_{4}\\
\frac{\partial\xi_{3}}{\partial u} & = & \ -s\xi_{2}+\e^{r}\sin(\theta)\xi_{1}-\e^{-r}\cos(\theta)\xi_{4}\\
\frac{\partial\xi_{4}}{\partial u} & = & \ \e^{r}(-\sin(\theta)\xi_{5}-\cos(\theta)\xi_{6})+\e^{-r}(\cos(\theta)\xi_{3}-\sin(\theta)\xi_{2})\\
\frac{\partial\xi_{5}}{\partial u} & = & \ s\xi_{6}+\e^{r}\sin(\theta)\xi_{4}-\e^{-r}\cos(\theta)\xi_{1}\\
\frac{\partial\xi_{6}}{\partial u} & = & \ -s\xi_{5}+\e^{r}\cos(\theta)\xi_{4}-\e^{-r}\sin(\theta)\xi_{1}\\
\frac{\partial r}{\partial u} & = & \ \xi_{1}\\
\frac{\partial s}{\partial u} & = & \ \e^{r}(-\cos(\theta)\xi_{3}-\sin(\theta)\xi_{2})+\e^{-r}(\cos(\theta)\xi_{6}-\sin(\theta)\xi_{5})\end{eqnarray*}
\begin{eqnarray*}
\frac{\partial\tilde{\xi}_{1}}{\partial u} & = & \ \e^{r}(\cos(\theta)\tilde{\xi}_{2}-\sin(\theta)\tilde{\xi}_{3})+\e^{-r}(\cos(\theta)\tilde{\xi}_{5}+\sin(\theta)\tilde{\xi}_{6})\\
\frac{\partial\tilde{\xi}_{2}}{\partial u} & = & \ s\tilde{\xi}_{3}-\e^{r}\cos(\theta)\tilde{\xi}_{1}+\e^{-r}\sin(\theta)\tilde{\xi}_{4}\\
\frac{\partial\tilde{\xi}_{3}}{\partial u} & = & \ -s\tilde{\xi}_{2}+\e^{r}\sin(\theta)\tilde{\xi}_{1}-\e^{-r}\cos(\theta)\tilde{\xi}_{4}\\
\frac{\partial\tilde{\xi}_{4}}{\partial u} & = & \ \e^{r}(-\sin(\theta)\tilde{\xi}_{5}-\cos(\theta)\tilde{\xi}_{6})+\e^{-r}(\cos(\theta)\tilde{\xi}_{3}-\sin(\theta)\tilde{\xi}_{2})\\
\frac{\partial\tilde{\xi}_{5}}{\partial u} & = & \ s\tilde{\xi}_{6}+\e^{r}\sin(\theta)\tilde{\xi}_{4}-\e^{-r}\cos(\theta)\tilde{\xi}_{1}\\
\frac{\partial\tilde{\xi}_{6}}{\partial u} & = & \ -s\tilde{\xi}_{5}+\e^{r}\cos(\theta)\tilde{\xi}_{4}-\e^{-r}\sin(\theta)\tilde{\xi}_{1}\end{eqnarray*}
and,\begin{eqnarray}
\frac{\partial\xi_{1}}{\partial v} & = & \ -\e^{r}(\cos(\theta)\xi_{3}+\sin(\theta)\xi_{2})+\e^{-r}(\cos(\theta)\xi_{6}-\sin(\theta)\xi_{5})\nonumber \\
\frac{\partial\xi_{2}}{\partial v} & = & \ -\xi_{1}\xi_{3}+\e^{r}\sin(\theta)\xi_{1}+\e^{-r}\cos(\theta)\xi_{4}\nonumber \\
\frac{\partial\xi_{3}}{\partial v} & = & \ \xi_{1}\xi_{2}+\e^{r}\cos(\theta)\xi_{1}+\e^{-r}\sin(\theta)\xi_{4}\nonumber \\
\frac{\partial\xi_{4}}{\partial v} & = & \ \e^{r}(\sin(\theta)\xi_{6}-\cos(\theta)\xi_{5})+\e^{-r}(-\cos(\theta)\xi_{2}-\sin(\theta)\xi_{3})\nonumber \\
\frac{\partial\xi_{5}}{\partial v} & = & \ -\xi_{1}\xi_{6}+\e^{r}\cos(\theta)\xi_{4}+\e^{-r}\sin(\theta)\xi_{1}\nonumber \\
\frac{\partial\xi_{6}}{\partial v} & = & \ \xi_{1}\xi_{5}-\e^{r}\sin(\theta)\xi_{4}-\e^{-r}\cos(\theta)\xi_{1}\nonumber \\
\frac{\partial r}{\partial v} & = & \ s\label{eq:2}\\
\frac{\partial s}{\partial v} & = & \ -2\sinh(2r)+\e^{r}(\sin(\theta)\xi_{3}-\cos(\theta)\xi_{2})+\e^{-r}(-\sin(\theta)\xi_{6}-\cos(\theta)\xi_{5})\nonumber \end{eqnarray}
\begin{eqnarray*}
\frac{\partial\tilde{\xi}_{1}}{\partial v} & = & \ -\e^{r}(\cos(\theta)\tilde{\xi}_{3}+\sin(\theta)\tilde{\xi}_{2})+\e^{-r}(\cos(\theta)\tilde{\xi}_{6}-\sin(\theta)\tilde{\xi}_{5})\\
\frac{\partial\tilde{\xi}_{2}}{\partial v} & = & \ -\xi_{1}\tilde{\xi}_{3}+\e^{r}\sin(\theta)\tilde{\xi}_{1}+\e^{-r}\cos(\theta)\tilde{\xi}_{4}\\
\frac{\partial\tilde{\xi}_{3}}{\partial v} & = & \ \xi_{1}\tilde{\xi}_{2}+\e^{r}\cos(\theta)\tilde{\xi}_{1}+\e^{-r}\sin(\theta)\tilde{\xi}_{4}\\
\frac{\partial\tilde{\xi}_{4}}{\partial v} & = & \ \e^{r}(\sin(\theta)\tilde{\xi}_{6}-\cos(\theta)\tilde{\xi}_{5})+\e^{-r}(-\cos(\theta)\tilde{\xi}_{2}-\sin(\theta)\tilde{\xi}_{3})\\
\frac{\partial\tilde{\xi}_{5}}{\partial v} & = & \ -\xi_{1}\tilde{\xi}_{6}+\e^{r}\cos(\theta)\tilde{\xi}_{4}+\e^{-r}\sin(\theta)\tilde{\xi}_{1}\\
\frac{\partial\tilde{\xi}_{6}}{\partial v} & = & \ \xi_{1}\tilde{\xi}_{5}-\e^{r}\sin(\theta)\tilde{\xi}_{4}-\e^{-r}\cos(\theta)\tilde{\xi}_{1}\end{eqnarray*}
Moreover, The system given by the equations (\ref{eq:2}) is integrable.
\end{thm}
\begin{proof}
This is long direct computation. 
\end{proof}

\subsection{First integrals and existence of global solutions}

In this subsection we will prove that the solutions of the sinh-Gordon
equations given by the integrable system (\ref{eq:1}) are defined
in the whole of $\R^{2}$. In order to prove this, we first establish
some lemmas.

\begin{lem}
\label{lem:M}For a given solution of the system (\ref{eq:1}), the
functions $\xi_{1},\dots,\xi_{6}$ defined in Theorem (\ref{thm:system-2})
satisfy\[
M=\frac{1}{2}\{\xi_{1}^{2}+\cdots+\xi_{6}^{2}\}\]
is a constant.
\end{lem}
\begin{proof}
A direct computation using Theorem (\ref{thm:system-2}) gives us
that\begin{eqnarray*}
\frac{\partial M}{\partial u} & = & \xi_{1}\frac{\partial\xi_{1}}{\partial u}+\cdots+\xi_{6}\frac{\partial\xi_{6}}{\partial u}\\
 & = & \xi_{1}(\e^{r}(\cos(\theta)\xi_{2}-\sin(\theta)\xi_{3})+\e^{-r}(\cos(\theta)\xi_{5}+\sin(\theta)\xi_{6}))\\
 &  & +\xi_{2}(s\xi_{3}-\e^{r}\cos(\theta)\xi_{1}+\e^{-r}\sin(\theta)\xi_{4})\\
 &  & +\xi_{3}(-s\xi_{2}+\e^{r}\sin(\theta)\xi_{1}-\e^{-r}\cos(\theta)\xi_{4})\\
 &  & +\xi_{4}(\e^{r}(-\sin(\theta)\xi_{5}-\cos(\theta)\xi_{6})+\e^{-r}(\cos(\theta)\xi_{3}-\sin(\theta)\xi_{2}))\\
 &  & +\xi_{5}(s\xi_{6}+\e^{r}\sin(\theta)\xi_{4}-\e^{-r}\cos(\theta)\xi_{1})\\
 &  & +\xi_{6}(-s\xi_{5}+\e^{r}\cos(\theta)\xi_{4}-\e^{-r}\sin(\theta)\xi_{1})\\
 & = & 0.\end{eqnarray*}
Similarly, \begin{eqnarray*}
\frac{\partial M}{\partial v} & = & \xi_{1}\frac{\partial\xi_{1}}{\partial v}+\cdots+\xi_{6}\frac{\partial\xi_{6}}{\partial v}\\
 & = & \xi_{1}(-\e^{r}(\cos(\theta)\xi_{3}+\sin(\theta)\xi_{2})+\e^{-r}(\cos(\theta)\xi_{6}-\sin(\theta)\xi_{5}))\\
 &  & +\xi_{2}(-\xi_{1}\xi_{3}+\e^{r}\sin(\theta)\xi_{1}+\e^{-r}\cos(\theta)\xi_{4})\\
 &  & +\xi_{3}(\xi_{1}\xi_{2}+\e^{r}\cos(\theta)\xi_{1}+\e^{-r}\sin(\theta)\xi_{4})\\
 &  & +\xi_{4}(\e^{-r}(-\cos(\theta)\xi_{2}-\sin(\theta)\xi_{3}))\\
 &  & +\xi_{5}(-\xi_{1}\xi_{6}+\e^{r}\cos(\theta)\xi_{4}+\e^{-r}\sin(\theta)\xi_{1})\\
 &  & +\xi_{6}(\xi_{1}\xi_{5}-\e^{r}\sin(\theta)\xi_{4}-\e^{-r}\cos(\theta)\xi_{1})\\
 & = & 0,\end{eqnarray*}
therefore, $M$ is a constant. 
\end{proof}
\begin{lem}
\label{lem:E}For a given solution of the system (\ref{eq:1}), the
functions $p,V_{1},V_{2},\nu$ satisfy\[
E=\frac{1}{2}\{\<p,p\>+\<V_{1},V_{1}\>+\<V_{2},V_{2}\>+\<\nu,\nu\>\}\]
is a constant.
\end{lem}
\begin{proof}
As in the proof of the previous lemma, a direct computation shows
that $\frac{\partial E}{\partial u}=\frac{\partial E}{\partial v}=0$. 
\end{proof}
\begin{lem}
\label{lem:A}For a given solution of the system (\ref{eq:1}), the
functions $\xi_{1},\dots,\xi_{6}$ defined in Theorem (\ref{thm:system-2})
satisfy\[
A=\e^{r}(\cos(\theta)\xi_{2}-\sin(\theta)\xi_{3})-\e^{-r}(\cos(\theta)\xi_{5}+\sin(\theta)\xi_{6})+\frac{1}{2}s^{2}+\cosh(2r)-\frac{1}{2}\left(\xi_{1}\right)^{2}\]
is a constant.
\end{lem}
\begin{proof}
Similarly to the previous two lemmas, we prove that $\frac{\partial A}{\partial u}=\frac{\partial A}{\partial v}=0$.

Denote by \begin{eqnarray*}
B & = & \e^{r}(\cos(\theta)\xi_{2}-\sin(\theta)\xi_{3})-\e^{-r}(\cos(\theta)\xi_{5}+\sin(\theta)\xi_{6})\ \com{and}\\
C & = & \frac{\partial\xi_{1}}{\partial u}=\e^{r}(\cos(\theta)\xi_{2}-\sin(\theta)\xi_{3})+\e^{-r}(\cos(\theta)\xi_{5}+\sin(\theta)\xi_{6}).\end{eqnarray*}
Notice that $B+\frac{1}{2}s^{2}-\frac{1}{2}\xi_{1}^{2}+\cosh(2r)=A$.
A direct computation shows that \begin{eqnarray*}
\frac{\partial B}{\partial u} & = & \xi_{1}C+\e^{r}\{\cos(\theta)(s\xi_{3}-\e^{r}\cos(\theta)\xi_{1}+\e^{-r}\sin(\theta)\xi_{4})\\
 &  & -\sin(\theta)(-s\xi_{2}+\e^{r}\sin(\theta)\xi_{1}-\e^{-r}\cos(\theta)\xi_{4})\}\\
 &  & -\e^{-r}\{\cos(\theta)(s\xi_{6}+\e^{r}\sin(\theta)\xi_{4}-\e^{-r}\cos(\theta)\xi_{1})\\
 &  & +\sin(\theta)(-s\xi_{5}+\e^{r}\cos(\theta)\xi_{4}-\e^{-r}\sin(\theta)\xi_{1})\}\\
 & = & \xi_{1}\frac{\partial\xi_{1}}{\partial u}+s(\e^{r}\cos(\theta)\xi_{3}+\e^{r}\sin(\theta)\xi_{2}-\e^{-r}\cos(\theta)\xi_{6}+\e^{-r}\sin(\theta)\xi_{5})\\
 &  & +\xi_{4}(\cos(\theta)\sin(\theta)+\cos(\theta)\sin(\theta)-\cos(\theta)\sin(\theta)-\cos(\theta)\sin(\theta))\\
 &  & +\xi_{1}(-\e^{2r}\cos^{2}(\theta)\xi_{3}-\e^{2r}\sin^{2}(\theta)\xi_{2}+\e^{-2r}\cos^{2}(\theta)+\e^{-2r}\sin^{2}(\theta))\\
 & = & \xi_{1}\frac{\partial\xi_{1}}{\partial u}-s\frac{\partial s}{\partial u}-2\xi_{1}\sinh(2r)\\
 & = & \frac{1}{2}\frac{\partial\xi_{1}^{2}}{\partial u}-\frac{1}{2}\frac{\partial s^{2}}{\partial u}-\frac{\partial\cosh(2r)}{\partial u}.\end{eqnarray*}
Therefore $\frac{\partial A}{\partial u}=0$. Similarly,  \begin{eqnarray*}
\frac{\partial B}{\partial v} & = & sC+\e^{r}\{\cos(\theta)(-\xi_{1}\xi_{3}+\e^{r}\sin(\theta)\xi_{1}+\e^{-r}\cos(\theta)\xi_{4})\\
 &  & -\sin(\theta)(-\xi_{1}\xi_{2}+\e^{r}\cos(\theta)\xi_{1}+\e^{-r}\sin(\theta)\xi_{4})\}\\
 &  & -\e^{-r}\{\cos(\theta)(-\xi_{1}\xi_{6}+\e^{r}\cos(\theta)\xi_{4}+\e^{-r}\sin(\theta)\xi_{1})\\
 &  & +\sin(\theta)(\xi_{1}\xi_{5}-\e^{r}\sin(\theta)\xi_{4}-\e^{-r}\cos(\theta)\xi_{1})\}\\
 & = & s(-2\sinh(2r)-\frac{\partial s}{\partial v})+\xi_{1}(-\e^{r}\cos(\theta)\xi_{3}+\e^{2r}\cos(\theta)\sin(\theta)\\
 &  & -\e^{2r}\sin(\theta)\cos(\theta)-\e^{r}\sin(\theta)\xi_{2})\\
 &  & +\e^{-r}\cos(\theta)\xi_{6}-\e^{-2r}\cos(\theta)\sin(\theta)-\e^{-r}\sin(\theta)\xi_{5}+\e^{-2r}\sin(\theta)\cos(\theta)\\
 &  & +\xi_{4}(\cos^{2}(\theta)-\sin^{2}(\theta)+\cos^{2}(\theta)+\sin^{2}(\theta))\\
 & = & -\frac{1}{2}\frac{\partial s^{2}}{\partial v}-\frac{\partial\cosh(2r)}{\partial v}+\frac{1}{2}\frac{\partial\xi_{1}^{2}}{\partial v}.\end{eqnarray*}

\end{proof}
\begin{cor}
\label{cor:R}Given a solution of the system (\ref{eq:1}). If $M$
and $A$ are the constants given by Lemmas (\ref{lem:M}) and (\ref{lem:A}),
respectively, if $(u_{0},v_{0})$ is any point in the domain of the
solution, and if $R$ is a real number such that\[
\cosh(2R)>A+4M\cosh(R)+\frac{M^{2}}{2}\quad\com{and}\quad R>|r(u_{0},v_{0})|\]
Then, $|r(u,v)|<R$ and\[
\frac{1}{2}s^{2}(u,v)+\cosh(2r(u,v))\le A+\frac{M^{2}}{2}+\cosh(2R)\]
for any $(u,v)$ in the domain of the solution.
\end{cor}
\begin{proof}
We have that\begin{eqnarray*}
\frac{1}{2}s^{2}(u,v)+\cosh(2r(u,v)) & = & A+\frac{1}{2}\xi_{1}^{2}+\e^{-r}(\cos(\theta)\xi_{5}+\sin(\theta)\xi_{6})-\e^{r}(\cos(\theta)\xi_{2}-\sin(\theta)\xi_{3})\\
 & \le & A+\frac{M^{2}}{2}+4M\cosh(r).\end{eqnarray*}
This inequality above shows that the result will follow once we prove
that $|r(u,v)|\le R$. We prove that $|r(u,v)|<R$ by contradiction.
If, for some $(u,v)$, $|r(u,v)|=R$, then, the inequality above implies
that at that $(u,v)$,\[
\cosh(2R)\le A+\frac{M^{2}}{2}+4M\cosh(R).\]
This is a contradiction with the choice of $R$ given in the hypotheses.
Therefore the Corollary follows.
\end{proof}
\begin{thm}
\label{thm:globally-defined}Any solution of the system (\ref{eq:1})
is defined on the entirety of $\R^{2}$.
\end{thm}
\begin{proof}
By Lemma (\ref{lem:M}), Lemma (\ref{lem:E}), and Corollary (\ref{cor:R}),
the solution of the system (\ref{eq:1}) remains bounded in $\R^{{18}}$
for all $(u,v)$, guaranteeing the existence of the solution for all
$(u,v)$. 
\end{proof}

\subsection{The integrable system and minimal immersion of the plane.}

In this section we prove that if we choose ${\bf x_{0}}$ properly,
a solution of the system (\ref{eq:1}) produces an example of a minimal
immersion of $\R^{2}$ into $S^{3}$. We start by showing that the
vector fields $Z$ and $W$ that define the system (\ref{eq:1}) define
vector fields in the manifold $SO(4)\times\R^{2}$. To see this, consider
the set $SO(4)\times\R^{2}$ as the following subset of $\R^{{18}}$:\[
SO(4)\times\R^{2}=\{(V_{0},V_{1},V_{2},V_{3},(r,s))\in\left(\R^{4}\right)^{4}\times\R^{2}:\<V_{i},V_{j}\>=\delta_{ij},\, det\left[V_{0}|\cdots|V_{3}\right]=1\}.\]

With this in mind, it is not difficult to verify the following lemma.

\begin{lem}
\label{lem:Y}If a vector field $Y:\R^{{18}}\to\R^{{18}}$ can be
written as\begin{eqnarray}
Y(x) & = & \ (c_{12}(x)V_{1}+c_{13}(x)V_{2}+c_{14}(x)V_{3},-c_{12}(x)V_{0}+c_{23}(x)V_{2}+c_{24}(x)V_{4},\nonumber \\
 &  & \ -c_{13}(x)V_{0}-c_{23}(x)V_{1}+c_{34}(x)V_{3},-c_{14}(x)V_{0}-c_{24}(x)V_{1}-c_{34}(x)V_{2},\nonumber \\
 &  & \ f(x),\ g(x))\ \label{eq:Y}\\
 & = & (C(x)[V_{0}^{T},V_{1}^{T},V_{2}^{T},V_{3}^{T}]\ ,f(x),g(x)),\nonumber \end{eqnarray}
where $x=(V_{0},V_{1},V_{2},V_{3},(r,s))$ denotes a point in $\R^{{18}}$,
and $C:\R^{{18}}\to so(4)$, $f:\R^{{18}}\to\R$ and $g:\R^{{18}}\to\R$
are smooth functions, then the restriction of $Y$ to $SO(4)\times\R^{2}$
defines a vector field tangent to $SO(4)\times\R^{2}$.
\end{lem}
As a consequence of this lemma we have:

\begin{lem}
The vector fields $Z$ and $W$ that define the integrable system
(\ref{eq:1}) define tangent vector fields in $SO(4)\times\R^{2}$.
\end{lem}
\begin{proof}
Notice that the notation used in Lemma (\ref{lem:Y}) and the notation
used in the system (\ref{eq:1}) are the same after identifying $p=V_{0}$
and $\nu=V_{3}$. The lemma follows because the vector field $Z$
can be written in the form of Lemma (\ref{lem:Y}) with\[
C_{1}(x)=\left(\begin{array}{cccc}
0 & e^{-r}\cos(\theta) & \e^{-r}\sin(\theta) & 0\\
-\e^{-r}\cos(\theta) & 0 & s & \e^{r}\cos(\theta)\\
-\e^{-r}\sin(\theta) & -s & 0 & -\e^{r}\sin(\theta)\\
0 & -\e^{r}\cos(\theta) & \e^{r}\sin(\theta) & 0\end{array}\right)\]
and\begin{eqnarray*}
f_{1}(p,V_{1},V_{2},\nu,(r,s)) & = & \<Bp,\nu\>\quad\quad\quad\com{and}\\
g_{1}(p,V_{1},V_{2},\nu,(r,s)) & = & \cos(\theta)\<BV_{2},\e^{-r}\nu-\e^{r}p\>-\sin(\theta)\<BV_{1},\e^{-r}\nu+\e^{r}p\>.\end{eqnarray*}
Also, the vector field $W$ can be written in the form of the lemma
(\ref{lem:Y}) with\[
C_{2}(x)=\left(\begin{array}{cccc}
0 & -\e^{-r}\sin(\theta) & \e^{-r}\cos(\theta) & 0\\
\e^{-r}\sin(\theta) & 0 & -\<Bp,\nu\> & -\e^{r}\sin(\theta)\\
-\e^{-r}\cos(\theta) & \<Bp,\nu\> & 0 & \e^{r}\cos(\theta)\cr0\\
0 & \e^{r}\sin(\theta) & \e^{r}\cos(\theta) & 0\end{array}\right)\]
and\begin{eqnarray*}
f_{2}(p,V_{1},V_{2},\nu,(r,s)) & = & s\quad\quad\quad\com{and}\\
g_{2}(p,V_{1},V_{2},\nu,(r,s)) & = & -2\sinh(2r)-\sin(\theta)\<BV_{2},\e^{-r}\nu-\e^{r}p\>-\cos(\theta)\<BV_{1},\e^{-r}\nu+\e^{r}p\>.\end{eqnarray*}
Lemma (\ref{lem:Y}) then implies that $Z$ and $W$ are tangent to
$SO(4)\times\R^{2}$.
\end{proof}
\begin{thm}
\label{thm:solutions-of-system}Let, $Z$ and $W$ be the vector fields
that defined the integrable system (\ref{eq:1}) and let ${\bf x^{0}}=(p^{0},V_{1}^{0},V_{1}^{0},\nu^{0},r^{0},s^{0})\in\R^{{18}}$
be such that $\{p^{0},V_{1}^{0},V_{1}^{0},\nu^{0}\}$ is an orthonormal
basis of $\R^{4}$ and let\[
\Theta_{Z}:\R\times\R^{{18}}\to\R^{{18}}\quad\hbox{and}\quad\Theta_{W}:\R\times\R^{{18}}\to\R^{{18}}\]
be the flows of the vector fields $Z$ and $W$ respectively. If $\phi:\R^{2}\to\R^{{18}}$
is given by\[
\phi(u,v)=\Theta_{Z}(u,\Theta_{W}(v,{\bf x^{0}}))=(\phi_{1}(u,v),\dots,\phi_{18}(u,v))\]
then the map\[
\rho:\R^{2}\to\R^{4}\ \hbox{given\, by}\ \rho(u,v)=(\phi_{1}(u,v),\dots,\phi_{4}(u,v))\]
satisfies $M=\rho(\R^{2})\subset S^{3}$, $\rho$ is a minimal immersion
of $S^{3}$ with principal curvature at $\rho(u,v)$ given by $\pm\e^{2\phi_{17}(u,v)}$.
More precisely, the map\[
\nu:\R^{2}\to\R^{4}\ \hbox{given\, by}\ \nu(u,v)=(\phi_{13}(u,v),\phi_{14}(u,v),\phi_{15}(u,v),\phi_{16}(u,v))\]
is the Gauss map of the immersion $\rho$, and\[
(\phi_{5}(u,v),\phi_{6}(u,v),\phi_{7}(u,v),\phi_{8}(u,v))\ \hbox{and}\ (\phi_{9}(u,v),\phi_{10}(u,v),\phi_{11}(u,v),\phi_{12}(u,v))\]
are the principal directions of the immersion $\rho$.
\end{thm}
\begin{proof}
Denote this solution by $\phi=(p,V_{1},V_{2},\nu,r,s)$, where $p=\rho,V_{1},V_{2},\nu:\R^{2}\to\R^{4}$
and $r,s:\R^{2}\to\R$. By the previous lemma, $|p|=|\rho|=1$, because
the initial conditions belong to $SO(4)\times\R^{2}$ and therefore
the whole solution stays in $SO(4)\times\R^{2}$. By the form of the
vector field $Z$ and $W$ or, equivalently, by the fact that $\phi$
is a solution of the system (\ref{eq:1}), \[
\frac{\partial p}{\partial u}=\e^{-r}(\cos(\theta)V_{1}+\sin(\theta)V_{2})\quad\com{and}\quad\frac{\partial p}{\partial v}=\e^{-r}(\cos(\theta)V_{2}-\sin(\theta)V_{1}).\]
Again by the fact that the solution remains in $SO(4)\times\R^{2}$,
the first fundamental form of the parameterized surface $\rho=p:\R^{2}\to S^{3}$
is given by\[
E=\<\frac{\partial p}{\partial u},\frac{\partial p}{\partial u}\>=\e^{-2r},\ F=\<\frac{\partial p}{\partial u},\frac{\partial p}{\partial v}\>=0,\ \com{and}\ G=\<\frac{\partial p}{\partial u},\frac{\partial p}{\partial v}\>=\e^{-2r}.\]
Therefore the immersion $\rho$ is a conformal immersion and for every
$(u,v)\in\R^{2}$, the vectors $\left\{ V_{1}(u,v),V_{2}(u,v)\right\} $
form a basis of the tangent space $T_{\rho(u,v)}M$. More precisely,
\[
V_{1}=\e^{r}(\cos(\theta)\frac{\partial p}{\partial u}-\sin(\theta)\frac{\partial p}{\partial v})\quad\com{and}\quad V_{2}=\e^{r}(\sin(\theta)\frac{\partial p}{\partial u}+\cos(\theta)\frac{\partial p}{\partial v}).\]
Once again from the fact that the solution remains in $SO(4)\times\R^{2}$,
 the map $\nu:\R^{2}\to S^{3}$ defines the Gauss map of the immersion
$\rho$. Since $\phi$ is a solution of the system (\ref{eq:1}),\[
\frac{\partial\nu}{\partial u}=d\nu(\frac{\partial p}{\partial u})=\e^{r}(-\cos(\theta)V_{1}+\sin(\theta)V_{2})\quad\com{and}\quad\frac{\partial\nu}{\partial v}=d\nu(\frac{\partial p}{\partial v})=\e^{r}(\sin(\theta)V_{1}+\cos(\theta)V_{2}).\]
In the previous equalities we identify the Gauss map defined in $\R^{2}$
with the Gauss map defined in $M=\rho(\R^{2})\subset S^{3}$. The
previous equation implies that\begin{eqnarray*}
d\nu(V_{1}) & = & d\nu(\e^{r}(\cos(\theta)\frac{\partial p}{\partial u}-\sin(\theta)\frac{\partial p}{\partial v}))\\
 & = & \e^{r}\{\cos(\theta)(\e^{r}(-\cos(\theta)V_{1}+\sin(\theta)V_{2}))-\sin(\theta)(\e^{r}(\sin(\theta)V_{1}+\cos(\theta)V_{2}))\}\\
 & = & -\e^{2r}V_{1}.\end{eqnarray*}
Similarly, \begin{eqnarray*}
d\nu(V_{2}) & = & d\nu(\e^{r}(\sin(\theta)\frac{\partial p}{\partial u}+\cos(\theta)\frac{\partial p}{\partial v}))\\
 & = & \e^{r}\{\sin(\theta)(\e^{r}(-\cos(\theta)V_{1}+\sin(\theta)V_{2}))+\cos(\theta)(\e^{r}(\sin(\theta)V_{1}+\cos(\theta)V_{2}))\}\\
 & = & e^{2r}V_{2}.\end{eqnarray*}
The previous two equalities show that the vectors $V_{1}$ and $V_{2}$
define principal directions and that the principal curvatures of the
immersion at the point $\rho(u,v)$ are $\pm e^{2r(u,v)}$. This completes
the proof of the Theorem.
\end{proof}
We also have that certain minimal immersions of tori induce solutions
of the system (1). The following theorem shows exactly which minimal
tori in $S^{3}$ are characterize by the integrable system (1).

\begin{thm}
\label{thm:nnt<8}Let $\tilde{\rho}:M\to S^{3}$ be a minimal immersed
torus in $S^{3}$. Using the notation given in the introduction, if
for some angle $\theta$ and some matrix $B\in so(4)$, $h_{\theta}=2f_{B}$,
then, it is possible to choose a covering map $\tau:\R^{2}\to M$,
maps $\rho:\R^{2}\to S^{3}$, $\nu:\R^{2}\to S^{3}$, $V_{1},V_{2}:\R^{2}\to S^{3}$,
and a function $r:\R^{2}\to\R$ using Theorem (\ref{uniformization})
and its corollaries, such that\[
\phi(u,v)=(\rho(u,v),V_{1}(u,v),V_{2}(u,v),\nu(u,v),r(u,v),\frac{\partial r}{\partial v}(u,v))\]
is a solution of the system (\ref{eq:1}) with matrix $B$ and angle
$\theta$.
\end{thm}
\begin{proof}
Using Remark (\ref{rem:rotatecoord}), we can rotate coordinates so
that the maps $\rho$, $\nu$, $V_{1}$, and $V_{2}$ in Theorem (\ref{uniformization})
and Corollaries (\ref{cor:prin-dirs}) and (\ref{cor:a}) satisfy\[
V_{1}(u,v)=W_{1}(\tau(u,v)),\quad V_{2}(u,v)=W_{2}(\tau(u,v)),\quad\nu(u,v)=N(\tau(u,v))\quad\com{and}\quad\alpha=\theta,\]
with $a(\tau(u,v))=\e^{2r}$. Since $\alpha=\theta$, \[
V_{1}=\e^{r}(\cos(\theta)\frac{\partial\rho}{\partial u}-\sin(\theta)\frac{\partial\rho}{\partial v})\quad\com{and}\quad V_{2}=\e^{r}(\sin(\theta)\frac{\partial\rho}{\partial u}+\cos(\theta)\frac{\partial\rho}{\partial v}),\]
if $2f_{B}=h_{\theta}$, then\begin{eqnarray*}
2\<B\rho,\nu\> & = & \cos(\theta)\e^{-3r}(\e^{r}(\cos(\theta)\frac{\partial\rho}{\partial u}-\sin(\theta)\frac{\partial\rho}{\partial v}))(\e^{2r})\\
 &  & +\sin(\theta)\e^{-3r}(\e^{r}(\sin(\theta)\frac{\partial\rho}{\partial u}+\cos(\theta)\frac{\partial\rho}{\partial v}))(\e^{2r})\\
 & = & 2\frac{\partial r}{\partial u}\end{eqnarray*}
so that\begin{equation}
2\<B\rho,\nu\>=2\frac{\partial r}{\partial u}=h_{\theta}\label{eq:dr/du=h-theta}\end{equation}

and, similarly,\begin{equation}
2\frac{\partial r}{\partial v}=2s=h_{\theta+\frac{\pi}{2}}.\label{eq:dr/dv=h-theta+2pi}\end{equation}
From the formulas for $V_{1}$ and $V_{2}$ in Corollary (\ref{cor:prin-dirs}),
we have that\[
\frac{\partial\rho}{\partial u}=\e^{-r}(V_{1}\cos(\theta)+\sin(\theta)V_{2})\quad\com{and}\quad\frac{\partial\rho}{\partial v}=\e^{-r}(-V_{1}\sin(\theta)+\sin(\theta)V_{2}),\]
which verify the equations in the integrable system (\ref{eq:1}).
Also, using the equation above and the formula for $\frac{\partial\nu}{\partial u}$
and $\frac{\partial\nu}{\partial v}$ in Theorem (\ref{uniformization}),
we get that\[
\frac{\partial\nu}{\partial u}=\e^{r}(-\cos(\theta)V_{1}+\sin(\theta)V_{2})\quad\com{and}\quad\frac{\partial\nu}{\partial v}=\e^{r}(\sin(\theta)V_{1}+\cos(\theta)V_{2})\]
Which verify the equations in the integrable system (\ref{eq:1}).
In the same way a direct computation shows that derivatives of $\frac{\partial V_{i}}{\partial u}$
satisfy the equations of the the system (\ref{eq:1}). In order to
complete the proof of this lemma, let us check the equation for $\frac{\partial s}{\partial v}$.
We have that\begin{eqnarray*}
\frac{\partial s}{\partial v} & = & \frac{\partial^{2}r}{\partial v^{2}}=-2\sinh(2r)-\frac{\partial^{2}r}{\partial u^{2}}\\
 & = & \ -2\sinh(2r)-\frac{\partial}{\partial u}\<B\rho,\nu\>\\
 & = & \ -2\sinh(2r)-\<B\frac{\partial\rho}{\partial u},\nu\>-\<B\rho,\frac{\partial\nu}{\partial u}\>\\
 & = & \ -\sin(\theta)\<BV_{2},\e^{-r}\nu-\e^{r}p\>-\cos(\theta)\<BV_{1},\e^{-r}\nu+\e^{r}p\>,\end{eqnarray*}
which verifies the equation in the integrable system (\ref{eq:1}).
The equation for $\frac{\partial s}{\partial u}$ is similar. 
\end{proof}
\begin{rem}
\label{rem:h}Arguing in the same way we did in the proof of the previous
theorem we have that if

\[
\phi(u,v)=(\rho(u,v),V_{1}(u,v),V_{2}(u,v),\nu(u,v),r(u,v),s(u,v))\]
is a doubly-periodic solution of the integral system (1) and $M$
is the torus $\frac{\R^{2}}{\sim}$, then,

\[
h_{\theta}([(u,v)])=2\frac{\partial r}{\partial u}(u,v)\quad\com{and}\quad h_{\theta+\frac{\pi}{2}}([(u,v)])=2\frac{\partial r}{\partial v}(u,v)=2s.\]
Moreover, for any $4\times4$ skew-symmetric matrix $\tilde{B}$,
$f_{\tilde{B}}([(u,v)])=\<\tilde{B}\rho(u,v),\nu(u,v)\>$. Also, since
$\phi$ satisfies the integrable system (1), then $h_{\theta}=2f_{B}$.
\end{rem}

\section{The Lawson-Hsiang examples}

The Lawson-Hsiang tori examples are characterized as those immersed
minimal tori in $S^{3}$ that are preserved by a 1-parameter group
of ambient isometries \cite{LH}. This section will show that these
examples can be seen, first, as those immersed minimal tori for which
there exists a nonzero matrix $B\in so(4)$ such that the function
$f_{B}:M\to\R$ is identically zero. Then we show that all these examples
are included in our new construction. We show that these examples
define solutions of the integrable system (\ref{eq:1}) with data
matrix $B\in so(4)$ identically zero. Then, we will prove that if
a solution of the integrable system (\ref{eq:1}) with $B={\bf 0}$
defines a minimal torus, then this torus must be one of the examples
of Lawson and Hsiang.

\begin{prop}
\label{pro:f_B-vanish}If $\widetilde{\rho}:M\to S^{3}$ is an immersed
closed minimal surface, such that $f_{B}:M\to{\bf R}$ vanishes for
some $B\ne{\bf 0}$, then $\tilde{\rho}(M)$ is invariant under the
group $\{\e^{tB}:t\in{\bf R}\}$, so that $M$ is one of the examples
of Hsiang-Lawson.
\end{prop}
\begin{proof}
Let $X:S^{3}\to\R^{4}$ be the tangent vector field on $S^{3}$ given
by $X(p)=Bp$. Since $0=f_{B}(m)=\<B\tilde{\rho}(m),N(m)\>$, then
$X$ induces a unit tangent vector field on $M$. Therefore the integrals
curves of the vector field $X$ that start in $\tilde{\rho}(M)$ remains
in $\tilde{\rho}(M)$, i.e. if $\tilde{\rho}(m)\in\tilde{\rho}(M)$
then $\e^{tB}\tilde{\rho}(m)\in\tilde{\rho}(M)$. 
\end{proof}
\begin{prop}
\label{pro:theta}Let $\tilde{\rho}:M\to S^{3}$ be minimal immersion
of a torus. If $f_{B}:M\to\R$ vanishes, then, for some angle $\theta$,
$h_{\theta}:M\to\R$ vanishes, and $\tilde{\rho}$ corresponds to
a solution of system (1) with $nnt(M)\leq6$.
\end{prop}
\begin{proof}
As in the previous proposition, the vector field $X(m)=B\tilde{\rho}(m)$
defines a tangent vector field on $M$. Since the function $a$ is
invariant under isometries and $X$ is a Killing vector field, then
the function $X(a)$ is identically zero. We will prove the proposition
by showing that for some fixed angle $\theta$ and some fixed real
number $\lambda$, $X=\lambda a^{-\frac{1}{2}}(\cos(\theta)W_{1}+\sin(\theta)W_{2})$.
Choose maps $\rho,\nu,V_{1},V_{2}:\R^{2}\to S^{3}$, a covering $\tau:\R^{2}\to M$
and a function $r:\R^{2}\to\R$ using Theorem (\ref{uniformization}),
and its corollaries, such that\[
W_{1}(\tau(u,v))=V_{1}(u,v),\quad W_{2}(\tau(u,v))=V_{2}(u,v)\quad\hbox{and}\quad N(\tau(u,v))=\nu(u,v).\]
With this special parameterization of this torus and having in mind
that $a(\tau(u,v))=\e^{2r(u,v)}$, we have that $\alpha=0$ and \[
V_{1}=\e^{r}\frac{\partial\rho}{\partial u},\quad V_{2}=\e^{r}\frac{\partial\rho}{\partial v},\quad W_{1}(a)(\tau(u,v))=2\e^{3r(u,v)}\frac{\partial r}{\partial u}(u,v)\quad\com{and}\quad W_{2}(a)(\tau(u,v))=2\e^{3r(u,v)}\frac{\partial r}{\partial v}(u,v).\]
Since $X$ is a tangent vector field, $X(\tau(u,v))=f(u,v)V_{1}(u,v)+g(u,v)V_{2}(u,v)$
for two doubly periodic smooth functions $f,g:\R^{2}\to{\bf \R}$.
Since, moreover, $X$ is a Killing vector field,\[
\<\nabla_{W_{1}}X,W_{1}\>(\tau(u,v))=V_{1}(f)(u,v)-\frac{W_{2}(a)}{2a}(\tau(u,v))g(u,v)=\e^{r}(\frac{\partial f}{\partial u}-g\frac{\partial r}{\partial v})=0,\]
\[
\<\nabla_{W_{2}}X,W_{2}\>(\tau(u,v))=V_{2}(g)(u,v)-\frac{W_{1}(a)}{2a}(\tau(u,v))f(u,v)=\e^{r}(\frac{\partial g}{\partial v}-f\frac{\partial r}{\partial u})=0,\mbox{ and}\]
\begin{eqnarray*}
(\<\nabla_{W_{1}}X,W_{2}\>+\<\nabla_{W_{2}}X,W_{1}\>)(\tau(u,v)) & = & V_{1}(g)(u,v)+\frac{W_{2}(a)}{2a}(\tau(u,v))f(u,v)+\\
 &  & V_{2}(f)(u,v)+\frac{W_{1}(a)}{2a}(\tau(u,v))g(u,v)\\
 & = & \e^{r}(\frac{\partial g}{\partial u}+f\frac{\partial r}{\partial v}+\frac{\partial f}{\partial v}+g\frac{\partial r}{\partial u})\\
 & = & 0.\end{eqnarray*}
A direct verification gives that the three equations above imply that
the function $h(u+iv)=(\e^{r}f)(u,v)+i(\e^{r}g)(u,v)$ is an analytic
function. Since $h$ is doubly periodic in $\R^{2}$, and in particular
it is bounded, then we get that the function $h$ is constant. We
can write this constant as $\lambda\cos(\theta)+i\lambda\sin(\theta)$
with $\lambda\ne0$. 

The rest of the proposition follows from Theorem (\ref{thm:nnt<8}),
since in this case $nnt(M)\leq6$.
\end{proof}
The previous proposition shows that all the examples discussed in
\cite{LH} are included in the family given by the system (\ref{eq:1}).
More precisely, each one of them is included in one system with $B={\bf 0}$.
The following proposition shows that if $B={\bf 0}$, then every torus
in the system (\ref{eq:1}) is one of the examples in \cite{LH}.
Recall that our examples characterize those minimal immersions such
that $h_{\theta}=f_{B}$, therefore, the condition $B={\bf 0}$ implies
that $h_{\theta}$ vanishes for some angle $\theta$.

\begin{prop}
\label{pro:a-1-variable}Let $\tilde{\rho}:M\to S^{3}$ be a minimal
immersion of a torus. If for some $\theta$, $h_{\theta}:M\to\R$
vanishes, then $f_{B}$ vanishes for some nonzero skew-symmetric matrix
$B$.
\end{prop}
\begin{proof}
Define the vector field $X$ by $X=a^{-\frac{1}{2}}\cos(\theta)W_{1}+a^{-\frac{1}{2}}\sin(\theta)W_{2}$.
The following identities show that $X$ is a Killing vector field
on $M$.\begin{eqnarray*}
\<\nabla_{W_{1}}X,W_{1}\> & = & -\frac{1}{2}a^{-\frac{3}{2}}W_{1}(a)\cos(\theta)-a^{-\frac{1}{2}}\frac{1}{2a}W_{2}(a)\sin(\theta)=-\frac{1}{2a}h_{\theta}=0\\
\<\nabla_{W_{2}}X,W_{2}\> & = & -\frac{1}{2}a^{-\frac{3}{2}}W_{2}(a)\sin(\theta)-a^{-\frac{1}{2}}\frac{1}{2a}W_{1}(a)\cos(\theta)=-\frac{1}{2a}h_{\theta}=0\\
\<\nabla_{W_{1}}X,W_{2}\> & = & -\frac{1}{2}a^{-\frac{3}{2}}W_{1}(a)\sin(\theta)+a^{-\frac{1}{2}}\frac{1}{2a}W_{2}(a)\cos(\theta)\\
\<\nabla_{W_{2}}X,W_{1}\> & = & -\frac{1}{2}a^{-\frac{3}{2}}W_{2}(a)\cos(\theta)+a^{-\frac{1}{2}}\frac{1}{2a}W_{1}(a)\sin(\theta)=-\<\nabla_{W_{1}}X,W_{2}\>.\end{eqnarray*}
Therefore the map $\Theta_{X}(t,\ \cdot):M\to M$ defines a 1-parameter
group of isometries in $M$. By Theorem (\ref{thm:Ramanaham}), $M$
is invariant under a 1-parameter group of isometries of $S^{3}$,
and therefore $f_{B}$ vanishes for some nonzero $B\in so(4)$.
\end{proof}

\section{Minimal surfaces with natural nullity less than 8}

The examples of minimal tori found in \cite{LH} that are not Clifford
tori can be divided into three categories $F_{1},F_{2}$ and $F_{3}$.
The first one, $F_{1}$, consists of the immersions given by

\[
\tilde{\rho}(u,v)=(\cos(mx)\cos(y),\sin(mx)\cos(y),\cos(kx)\sin(y),\sin(kx)\sin(y))\]

Where $m$ and $k$ are two relatively-prime positive integers. These
examples can be characterized by the property that the principal curvature
function $a:M\to\R$ is constant along a direction that makes a constant
angle of $\frac{\pi}{4}$ with respect to one of the principal directions.
The second category, $F_{2}$, are the examples found initially by
Otsuki \cite{O}, and are characterized by the property that the function
$a$ is constant along one of the principal directions. The third
category, $F_{3}$, are the new examples found in the paper \cite{LH}
that complete the classification of minimal immersions of tori that
are invariant under a group of isometries of $S^{3}$.

Since there is an explicit parameterization $\tilde{\rho}$, these
examples explicitly give solutions for the system (\ref{eq:1}). 

\begin{prop}
\label{pro:HL-are-sols} Assume that the variables $x$ and $y$ are
related to the variables $u$ and $v$ by the following equations:\[
u=\int_{0}^{y}\sqrt{\frac{mk}{m^{2}\cos^{2}(t)+k^{2}\sin^{2}(t)}}dt\quad\com{and}\quad v=\sqrt{mk}x.\]
For any pair of positive real numbers $m$ and $k$, the map\[
\phi(u,v)=(\rho(u,v),V_{1}(u,v),V_{2}(u,v),\nu(u,v),r(u,v),s(u,v))\]
given by\begin{eqnarray*}
\rho(u,v) & = & (\cos(mx)\cos(y),\sin(mx)\cos(y),\cos(kx)\sin(y),\sin(kx)\sin(y))\\
V_{1}(u,v) & = & \frac{1}{\sqrt{2}}(-\cos(mx)\sin(y),-\sin(mx)\sin(y),\cos(kx)\cos(y),\sin(kx)\cos(y))\\
 &  & +\frac{1}{\sqrt{2(m^{2}\cos^{2}(y)+k^{2}\sin^{2}(y))}}(-m\sin(mx)\cos(y),m\cos(mx)\cos(y),\\
 &  & \qquad\qquad\qquad\qquad\qquad\qquad\qquad-k\sin(kx)\sin(y),k\cos(kx)\sin(y))\\
V_{2}(u,v) & = & \frac{1}{\sqrt{2}}(\cos(mx)\sin(y),\sin(mx)\sin(y),-\cos(kx)\cos(y),-\sin(kx)\cos(y))\\
 &  & +\frac{1}{\sqrt{2(m^{2}\cos^{2}(y)+k^{2}\sin^{2}(y))}}(-m\sin(mx)\cos(y),m\cos(mx)\cos(y),\\
 &  & \qquad\qquad\qquad\qquad\qquad\qquad\qquad-k\sin(kx)\sin(y),k\cos(kx)\sin(y))\\
\nu(u,v) & = & \sqrt{\frac{1}{m^{2}\cos^{2}(y)+k^{2}\sin^{2}(y)}}(k\sin(mx)\sin(y),-k\cos(mx)\sin(y),\\
 &  & \qquad\qquad\qquad\qquad\qquad\qquad-m\sin(kx)\cos(y),m\cos(kx)\cos(y))\end{eqnarray*}
and\[
r(u,v)=\frac{1}{2}\ln{\big{(}\frac{mk}{m^{2}\cos^{2}(y)+k^{2}\sin^{2}(y)}\big{)}},\quad s(u,v)=0\]
is a solution of the integrable system (\ref{eq:1}) with $\theta=-\frac{\pi}{4}$
and\[
B=\left(\begin{array}{cccc}
0 & \frac{m^{2}-k^{2}}{k\sqrt{mk}} & 0 & 0\\
-\frac{m^{2}-k^{2}}{k\sqrt{mk}} & 0 & 0 & 0\\
0 & 0 & 0 & 0\\
0 & 0 & 0 & 0\end{array}\right).\]

\end{prop}
Similarly, 

\begin{prop}
Assume that the variables $x$ and $y$ are related to the variables
$u$ and $v$ by the following equations:\[
v=\int_{0}^{y}\sqrt{\frac{mk}{m^{2}\cos^{2}(t)+k^{2}\sin^{2}(t)}}dt\quad\com{and}\quad u=\sqrt{mk}x\]
 For any pair of positive real numbers $m$ and $k$, the map\[
\phi(u,v)=(\rho(u,v),V_{1}(u,v),V_{2}(u,v),\nu(u,v),r(u,v),s(u,v))\]
given by\begin{eqnarray*}
\rho(u,v) & = & (\cos(mx)\cos(y),\sin(mx)\cos(y),\cos(kx)\sin(y),\sin(kx)\sin(y))\\
V_{1}(u,v) & = & \frac{1}{\sqrt{2}}(-\cos(mx)\sin(y),-\sin(mx)\sin(y),\cos(kx)\cos(y),\sin(kx)\cos(y))\\
 &  & +\frac{1}{\sqrt{2(m^{2}\cos^{2}(y)+k^{2}\sin^{2}(y))}}(-m\sin(mx)\cos(y),m\cos(mx)\cos(y),\\
 &  & \qquad\qquad\qquad\qquad\qquad\qquad-k\sin(kx)\sin(y),k\cos(kx)\sin(y))\\
V_{2}(u,v) & = & \frac{1}{\sqrt{2}}(\cos(mx)\sin(y),\sin(mx)\sin(y),-\cos(kx)\cos(y),-\sin(kx)\cos(y))\\
 &  & +\frac{1}{\sqrt{2(m^{2}\cos^{2}(y)+k^{2}\sin^{2}(y))}}(-m\sin(mx)\cos(y),m\cos(mx)\cos(y),\\
 &  & \qquad\qquad\qquad\qquad\qquad\qquad-k\sin(kx)\sin(y),k\cos(kx)\sin(y))\\
\nu(u,v) & = & \sqrt{\frac{1}{m^{2}\cos^{2}(y)+k^{2}\sin^{2}(y)}}(k\sin(mx)\sin(y),-k\cos(mx)\sin(y),\\
 &  & \qquad\qquad\qquad\qquad\qquad\qquad-m\sin(kx)\cos(y),m\cos(kx)\cos(y))\end{eqnarray*}
and\begin{eqnarray*}
r(u,v) & = & \frac{1}{2}\ln{\big{(}\frac{mk}{m^{2}\cos^{2}(y)+k^{2}\sin^{2}(y)}\big{)}},\\
s(u,v) & = & \frac{m^{2}-k^{2}}{\sqrt{mk(m^{2}\cos^{2}(y)+k^{2}\sin^{2}(y)})}\sin(y)\cos(y),\end{eqnarray*}
is a solution of the integrable system (\ref{eq:1}) with $\theta=-\frac{\pi}{4}$
and matrix $B={\bf 0}$.
\end{prop}
\begin{rem}
\label{rem:KS}From remark (\ref{rem:h}), for the examples we are
studying, those that come from solutions of the integrable system
(1), $h_{\theta}=2f_{B}$, therefore $h_{\theta}\in KS$. Since for
any $\theta$, $\hbox{span}\{h_{0},h_{\frac{\pi}{2}}\}=\hbox{span}\{h_{\theta},h_{\theta+\frac{\pi}{2}}\}$,
then, if $h_{\theta+\frac{\pi}{2}}=2f_{\tilde{B}}$ for some $4\times4$
skew-symmetric matrix $\tilde{B}$, we have that $\hbox{span}\{h_{0},h_{\frac{\pi}{2}}\}\subset KS$.
Notice that, again by remark (\ref{rem:h}), $h_{\theta+\frac{\pi}{2}}=2\frac{\partial r}{\partial v}(u,v)=s(u,v)$,
therefore, in these examples, showing that $\hbox{span}\{h_{0},h_{\frac{\pi}{2}}\}\subset KS$
is equivalent to showing that $s=\frac{\partial r}{\partial v}=f_{\tilde{B}}$
for some $4\times4$ skew-symmetric matrix $\tilde{B}$.
\end{rem}
From Proposition (\ref{pro:HL-are-sols}) we can deduce that the natural
nullity of the family $F_{1}$ is 5, since both of the functions $h_{0}$
and $h_{\frac{\pi}{2}}$ are contained in the space $\{f_{B}:B\in so(4)\}$.
Notice that by Proposition (\ref{pro:theta}), if $M$ is in one of
the families $F_{1}$, $F_{2}$, or $F_{3}$, then for some $\theta$
the function $h_{\theta}$ vanishes. The question that we will address
now is whether the function $h_{(\theta+\frac{\pi}{2})}$, in general,
is contained in the set $\{f_{B}:B\in so(4)\}$, as it is for tori
in $F_{1}$. Notice that the question is equivalent to whether the
natural nullity of the Lawson-Hsiang examples that are not Clifford
tori is equal to 5. The following two theorems resolve the question.

\renewcommand{\labelenumi}{(\alph{enumi})}

\begin{thm}
\label{thm:s-vanishes}Let $\phi:\R^{2}\to\R^{{18}}$ be a solution
of integrable system (\ref{eq:1}), and let $r(u,v)=\phi_{17}(u,v)$
and $s(u,v)=\phi_{18}(u,v)$. Assume that $\phi(0,0)=x^{0}=(e_{1},e_{2},e_{3},e_{4},r_{0},0)$
and $\frac{\partial r}{\partial u}(0,0)=0$. If \[
B=\left(\begin{array}{cccc}
0 & b_{1} & b_{2} & b_{3}\\
-b_{1} & 0 & b_{4} & b_{5}\\
-b_{2} & -b_{4} & 0 & b_{6}\\
-b_{3} & -b_{5} & -b_{6} & 0\end{array}\right),\]
then, $s$ vanishes if and only if $b_{3}=b_{4}=0$ and 
\begin{enumerate}
\item $-\e^{r_{0}}\cos(\theta)b_{1}+\e^{r_{0}}\sin(\theta)b_{2}+\e^{-r_{0}}\cos(\theta)b_{5}+\e^{-r_{0}}\sin(\theta)b_{6}=\ 2\sinh(2r_{0})$, 
\item $-\e^{r_{0}}\sin(\theta)b_{1}-\e^{r_{0}}\cos(\theta)b_{2}+\e^{-r_{0}}\sin(\theta)b_{5}-\e^{-r_{0}}\cos(\theta)b_{6}=0$,
and 
\item $-\e^{-r_{0}}\cos(\theta)b_{1}-\e^{-r_{0}}\sin(\theta)b_{2}+\e^{r_{0}}\cos(\theta)b_{5}-\e^{r_{0}}\sin(\theta)b_{6}=0$. 
\end{enumerate}
\end{thm}
\begin{proof}
We will use the integrable system (\ref{eq:2}) with $\tilde{B}={\bf 0}$.
Notice that\[
b_{3}=-\xi_{1}(0,0),\ b_{4}=-\xi_{4}(0,0),\ b_{1}=\xi_{2}(0,0),\ b_{2}=\xi_{3}(0,0),\ b_{5}=-\xi_{5}(0,0),\ b_{6}=-\xi_{6}(0,0).\]
Assume that $s(u,v)=0$ for every $(u,v)\in\R^{2}$. The equation
$b_{3}=0$ follows because we are assuming that $\frac{\partial r}{\partial u}(0,0)=\xi_{1}(0,0)=0$.
Equation (a) in the statement of the theorem follows from the equation
$\frac{\partial s}{\partial v}(0,0)=0$. Equation (b) follows from
the equation $\frac{\partial s}{\partial u}(0,0)=0$. We now prove
that $s\equiv0$ also implies that $b_{4}=0$ and equation (c) in
the statement of the theorem.

A direct computation shows the following two equations,\begin{eqnarray*}
\frac{\partial^{2}s}{\partial v\partial u} & = & \xi_{1}\big{(}-2\cosh(2r)+\e^{r}(\sin(\theta)\xi_{3}-\cos(\theta)\xi_{2})+\e^{-r}(\sin(\theta)\xi_{6}+\cos(\theta)\xi_{5})\big{)}\\
 &  & +s\big{(}-\e^{r}(\sin(\theta)\xi_{2}+\cos(\theta)\xi_{3})+\e^{-r}(\sin(\theta)\xi_{5}-\cos(\theta)\xi_{6})\big{)}-2\sin(2\theta)\xi_{4}\end{eqnarray*}
and\begin{eqnarray*}
\frac{\partial^{2}s}{\partial v^{2}} & = & s\big{(}-4\cosh(2r)+\e^{r}(\sin(\theta)\xi_{3}-\cos(\theta)\xi_{2})+\e^{-r}(\sin(\theta)\xi_{6}+\cos(\theta)\xi_{5})\big{)}\\
 &  & +\xi_{1}\big{(}\e^{r}(\sin(\theta)\xi_{2}+\cos(\theta)\xi_{3})+\e^{-r}(\cos(\theta)\xi_{6}-\sin(\theta)\xi_{5})\big{)}-2\cos(2\theta)\xi_{4}.\end{eqnarray*}
From the previous equations we get that $\xi_{4}(0,0)=-b_{4}=0$ and
that $\frac{\partial\xi_{4}}{\partial v}(0,0)=0\mbox{ because }\xi_{1}(0,0)=0$,
and\[
\frac{\partial\xi_{1}}{\partial v}(0,0)=\frac{\partial s}{\partial u}(0,0)=0.\]
A direct computation shows that the equation (c) in the statement
of the theorem is equivalent to the equation $\frac{\partial\xi_{4}}{\partial v}(0,0)=0$.
So we have shown one implication in the theorem.

We now show the other implication. Assume that we have the equations
(a), (b) and (c) on the statement of the theorem, and also $b_{4}=b_{3}=0$.
These 5 conditions are equivalent to the conditions\[
\xi_{1}(0,0)=0,\ \xi_{4}(0,0)=0,\ \frac{\partial\xi_{1}}{\partial v}(0,0)=\frac{\partial s}{\partial u}(0,0)=0,\frac{\partial s}{\partial v}(0,0)=0,\ \com{and}\ \frac{\partial\xi_{4}}{\partial v}(0,0)=0.\]
Notice also that by assumption we have that $s(0,0)=0$. Using the
integrable system (\ref{eq:2}) we can see that the initial conditions
above imply that\begin{equation}
\frac{\partial\xi_{i}}{\partial u}(0,0)=\frac{\partial\xi_{i}}{\partial v}(0,0)=0,\ \com{for}\, i=2,3,5,6,\label{eq:partials-of-xi}\end{equation}
and, also, we can prove by induction that given $n\ge1$, $k$ and
$l$ non-negative integers such that $k+l=n$, there exists a polynomial
$P=P(t_{1},\dots,t_{9})$ such that\[
\frac{\partial^{n}r}{\partial u^{l}\partial v^{k}}=P(\e^{r},\e^{-r},s,\xi_{1},\dots,\xi_{6}).\]
Along with the equations in (\ref{eq:partials-of-xi}), these equations
imply that\[
\frac{\partial^{n}s}{\partial u^{l}\partial v^{k+1}}(0,0)=\frac{\partial(\frac{\partial^{n}r}{\partial u^{l}\partial v^{k}})}{\partial v}(0,0)=\frac{\partial P(\e^{r},\e^{-r},s,\xi_{1},\dots,\xi_{6})}{\partial v}(0,0)=0.\]
In the last equation we have also used the hypothesis that $\frac{\partial\xi_{1}}{\partial v}(0,0)=\frac{\partial\xi_{4}}{\partial v}(0,0)=0$.
We should point out that we have used the fact that the function $r$
is real analytic, which follows from the fact that $\Delta r+2\sinh(2r)=0$.
\end{proof}
\begin{rem}
\label{rem:B^1,B^2}If $r_{0}=0$, then for any angle $\theta$, the
matrices $B$ that satisfy the conditions of the previous theorem
form a 2-dimensional subspace of $so(4)$. They have the following
form:\[
B_{b_{1},b_{2}}=\left(\begin{array}{cccc}
0 & b_{1} & b_{2} & 0\\
-b_{1} & 0 & 0 & b_{1}\\
-b_{2} & 0 & 0 & -b_{2}\\
0 & -b_{1} & b_{2} & 0\end{array}\right).\]
The reason for the existence of this two-dimensional subspace of $so(4)$
is that every Clifford torus $M$ that contains the point $p_{0}=(1,0,0,0)$
with tangent space containing the vectors $(0,1,0,0)$ and $(0,0,1,0)$,
i.e., with $N(p_{0})=(0,0,0,\pm1)$, have the property that $\e^{tB_{b_{1},b_{2}}}M=M$.
If $r_{0}\ne0$, then, given $\theta$, the matrices $B$ that satisfy
the conditions of the previous theorem form a 1-dimensional affine
space and are of the form \begin{eqnarray*}
B=B_{\theta}^{1}+\lambda B_{\theta}^{2} & = & \left(\begin{array}{cccc}
0 & -\e^{r_{0}}\cos(\theta) & \e^{r_{0}}\sin(\theta) & 0\\
\e^{r_{0}}\cos(\theta) & 0 & 0 & -\e^{-r_{0}}\cos(\theta)\\
-\e^{r_{0}}\sin(\theta) & 0 & 0 & -\e^{-r_{0}}\sin(\theta)\\
0 & \e^{-r_{0}}\cos(\theta) & \e^{-r_{0}}\sin(\theta) & 0\end{array}\right)\\
 &  & +\lambda\left(\begin{array}{cccc}
0 & \e^{-r_{0}}\sin(\theta) & -\e^{-r_{0}}\cos(\theta) & 0\\
-\e^{-r_{0}}\sin(\theta) & 0 & 0 & \e^{r_{0}}\sin(\theta)\\
\e^{-r_{0}}\cos(\theta) & 0 & 0 & \e^{r_{0}}\cos(\theta)\\
0 & -\e^{r_{0}}\sin(\theta) & -\e^{r_{0}}\cos(\theta) & 0\end{array}\right).\end{eqnarray*}

\end{rem}
\begin{thm}
\label{thm:nullity5}If $M\subset S^{3}$ is an immersed torus invariant
under a 1-parametric group of isometries of $S^{3}$, then $nnt(M)=kn(M)$
and therefore the natural nullity $nnt(M)\leq5$.
\end{thm}
\begin{proof}
By Proposition (\ref{pro:a-1-variable}) we know that for some angle
$\theta$, $(\cos(\theta)V_{1}+\sin(\theta)V_{2})(a)=0$ where $a:M\to{\bf R}$
is a positive function such that the principal curvatures of $M$
at $p$ are $\pm a(p)$. Without loss of generality, we will assume
that\[
e_{1}\in M,\ V_{1}(e_{1})=e_{2},\ V_{2}(e_{1})=e_{3},\nu(e_{1})=e_{2},\ \ln{a}(e_{1})=2r_{0},\ \com{and}\ \ \nabla a(e_{1})={\bf 0}.\]
Therefore, $M$ defines a solution of the system (\ref{eq:1}) associated
with the matrix $B={\bf 0}$ and $\theta$. Call this solution $\phi:\R^{2}\to\R^{{18}}$.
Without loss of generality we can assume that $\phi(0,0)=(e_{1},e_{2},e_{3},e_{4},r_{0},0)$.

Define $\tilde{\phi}$ to be the solution of the system (\ref{eq:1})
associated with a matrix $B=\{b_{ij}\}$ that satisfies the conditions
in the previous lemma and $\tilde{\theta}=\theta-\frac{\pi}{2}$.
Moreover we will take the initial solution that satisfies\[
\tilde{\phi}(0,0)=(e_{1},e_{2},e_{3},e_{4},r_{0},0).\]
Now consider the map $\hat{\phi}:\R^{2}\to\R^{{18}}$ given by

\begin{eqnarray*}
\hat{\phi}(u,v) & = & (\hat{\rho}(u,v),\hat{V_{1}}(u,v),\hat{V_{2}}(u,v),\hat{\nu}(u,v),\hat{r}(u,v),\hat{s}(u,v))\\
 & = & (\tilde{\rho}(-v,u),\tilde{V_{1}}(-v,u),\tilde{V_{2}}(-v,u),\tilde{\nu}(-v,u),\tilde{r}(-v,u),-\<B\tilde{\rho},\tilde{\nu}\>),\end{eqnarray*}
where \[
\tilde{\phi}(\tilde{u},\tilde{v})=(\tilde{\rho}(\tilde{u},\tilde{v}),\tilde{V_{1}}(\tilde{u},\tilde{v}),\tilde{V_{2}}(\tilde{u},\tilde{v}),\tilde{\nu}(\tilde{u},\tilde{v}),\tilde{r}(\tilde{u},\tilde{v}),\tilde{s}(\tilde{u},\tilde{v})).\]
It is clear that $\hat{\phi}(0,0)=(e_{1},e_{2},e_{3},e_{4},r_{0},0)$.
Notice that, by the way $B$ was chosen, we have that $\tilde{s}=0$
for every $(\tilde{u},\tilde{v})\in\R^{2}$. Also, a direct computation
shows that $\hat{\phi}$ is a solution of the system (\ref{eq:1})
with $B={\bf 0}$ and the angle $\theta$, therefore, $\hat{\phi}(u,v)=\phi(u,v)$,
and so\[
\frac{\partial r}{\partial v}=-\frac{\partial\tilde{r}}{\partial\tilde{u}}=-\<B\rho,\nu\>.\]
This equality is equivalent to the fact that $\sin(\theta)u_{1}-\cos(\theta)u_{2}=f_{B}$,
where the functions $u_{1}=h_{0},\: u_{2}=h_{\pi/2}$, and $f_{B}$
are defined in the first section. This last equation implies that
$h_{\theta+\frac{\pi}{2}}=-f_{B}$, therefore, $h_{\theta}$, which
is identically zero, and $h_{\theta+\frac{\pi}{2}}$ are functions
in $\{f_{C}:C\in so(4)\}$. Then, both functions $u_{1}$ and $u_{2}$
are also generated by the functions in the set $\{f_{C}:C\in so(4)\}$,
\emph{i.e.}, the natural nullity is $5$. Recall that the space $\{u_{C}:C\in so(4)\}$
is $5$-dimensional for any torus invariant under a 1-parameter group
of isometries in $S^{3}$. 
\end{proof}
\begin{cor}
\label{cor:nnt=00003D5-implies-ambient-isom}Let $M$ be torus with
natural nullity $5$, and let $\theta\in\R$ and $B\in so(4)$ be
such that $h_{\theta}=0$ and $h_{(\theta-\frac{\pi}{4})}=f_{B}$.
If $B_{\theta}^{1}$ and $B_{\theta}^{2}$ are defined as in the Remark
(\ref{rem:B^1,B^2}), then \[
B=B_{\theta}^{1}+\lambda_{0}B_{\theta}^{2}\quad\com{for\: some}\quad\lambda_{0}\quad\com{and}\quad\e^{\lambda B_{\theta}^{2}}M=M\quad\com{for\: every}\ \lambda\in\R.\]

\end{cor}
\begin{proof}
If $nnt(M)=5$, then certainly there is a nonzero $B\in so(4)$ so
that $f_{B}=0$. Proposition (\ref{pro:theta}) then implies the existence
of a $\theta$ for which $h_{\theta}=0$. $B=B_{\theta}^{1}+\lambda_{0}B_{\theta}^{2}$
follows from the previous theorem and the Remark (\ref{rem:B^1,B^2}).
The second part of the corollary follows from the fact that in the
argument used to prove Theorem (\ref{thm:nullity5}), we can choose
any $B$ that satisfies the conditions of Theorem (\ref{thm:s-vanishes}),
in particular if we can also choose $B=B_{\theta}^{1}+(\lambda_{0}+1)B_{\theta}^{2}$
we will get that\[
\frac{\partial r}{\partial v}=\<(B_{\theta}^{1}+\lambda_{0}B_{\theta}^{2})p,\nu\>=\<(B_{\theta}^{1}+(\lambda_{0}+1)B_{\theta}^{2})p,\nu\>.\]
This equation implies that $f_{B_{\theta}^{2}}=0$. The corollary
follows by the Proposition (\ref{pro:f_B-vanish}). 
\end{proof}
\begin{cor}
\label{cor:nnt=00003D5-iff-HL}If $M$ is a minimal immersed torus
in $S^{3}$, then $nnt(M)\leq5$ if and only if $M$ is one of the
examples of Hsiang and Lawson.
\end{cor}
\begin{proof}
If $M$ has $nnt(M)\leq5$, then $kn(M)\leq5$. Therefore, for some
nonzero skew-symmetric matrix $B$, $f_{B}$ vanishes. By Proposition(\ref{pro:f_B-vanish}),
$M$ will be invariant under a 1-parameter subgroup of the rigid motions
of $S^{3}$, which, following \cite{LH}, implies that $M$ is one
of Hsiang and Lawson's examples. On the other hand, since any of the
Hsiang-Lawson examples are preserved by a one-parameter subgroup of
$SO(4)$, there is a $B\in so(4)$ for which $f_{B}=0$. Then Theorem
(\ref{thm:nullity5}) implies $nnt(M)\leq5$.
\end{proof}
Theorem (\ref{thm:nullity5}) and Corollary (\ref{cor:nnt=00003D5-implies-ambient-isom})
address the question of the injectivity of the function that sends
any pair $(\theta,B,r_{0})$ to the minimal immersion of the plane
with initial conditions $(e_{1},e_{2},e_{3},e_{4},r_{0},0)$. The
following result is in the same direction.

\begin{prop}
If for some $\theta_{2}\ne\theta_{1}+n\pi$ for any integer $n$,
$h_{\theta_{1}}=f_{B_{1}}$ and $h_{\theta_{2}}=f_{B_{2}}$, then
the natural nullity of $M$ is less than 7.
\end{prop}
\begin{proof}
The equations in the Proposition implies that the space $\{\lambda h_{\theta}:\ \lambda,\theta\in\R\}$
is a subset of the space $\{f_{B}:B\in so(4)\}$ which has dimension
at most $6$. The proposition then follows.
\end{proof}
\begin{lem}
If a solution of (1) satisfies $r(0,0)=r_{0}$, $\xi_{1}(0,0)=s(0,0)=\xi_{4}(0,0)=0$,
then $r(u,v)=r(-u,-v)$.
\end{lem}
\begin{proof}
A direct computation using the system (\ref{eq:2}) shows that the
conditions $\xi_{1}(0,0)=s(0,0)=\xi_{4}(0,0)=0$ give\[
\frac{\partial\xi_{i}}{\partial u}(0,0)=\frac{\partial\xi_{i}}{\partial v}(0,0)=0\quad\com{for}\quad i=2,3,5,6.\]
Let $C^{\omega}(\R^{2})$ be the set of analytic functions on $\R^{2}$
and let $P_{0}$ be the ideal of $C^{\omega}(\R^{2})$ generated by
the functions $\{\e^{r},\e^{-r},\xi_{2},\xi_{3},\xi_{5},\xi_{6}\}$.
Given a nonnegative integer $k$, define $P_{k}$ as the set of functions
in $C^{\omega}(\R^{2})$ that can be written as a homogeneous polynomial
of degree $k$ in the variables $s,\:\xi_{1}$ and $\xi_{4}$ with
coefficients in $P_{0}$. A direct computation using the system (\ref{eq:2})
give us that if $f\in P_{0}$, then $\frac{\partial f}{\partial u}$
and $\frac{\partial f}{\partial v}$ are in $P_{1}$. In the same
way, if $f\in P_{k}$ then $\frac{\partial f}{\partial u}$ and $\frac{\partial f}{\partial v}$
are in $P_{k+1}+P_{k-1}$. Now with these observations in mind, we
proceed to show that the function $r$ satisfies $r(u,v)=r(-u,-v)$,
by showing that all the partial derivatives of odd order of the function
$r$ vanish at $(0,0)$. To achieve this we first notice that the
first derivatives of $r$, the functions $\xi_{1}$ and $s$ vanish
at $(0,0)$. Then, notice that the second derivatives of $r$, i.e.
the first derivatives of $s$ and $\xi_{1}$, are functions in $P_{0}$.
The last statement implies that the third derivatives of $r$ are
in $P_{1}$ and therefore they vanish at $(0,0)$. Once we know that
the third derivatives of $r$ are in $P_{1}$ we get that the fourth
derivatives or $r$ are in $P_{0}+P_{2}$. If we continue with this
process we notice that if $k$ is a positive even integer, then the
$k$-th derivatives of $r$ are functions in $P_{0}+P_{2}+\cdots+P_{k-2}$,
and in the case that $k$ is a odd integer greater that $1$, then,
the $k$-th derivatives of $r$ are in $P_{1}+P_{3}+\cdots+P_{k-2}$.
Now, since $\xi_{1}(0,0)=s(0,0)=\xi_{4}(0,0)=0$, the odd derivatives
of the function $r$ vanish at $(0,0)$. 
\end{proof}
\begin{thm}
\label{thm:nnt<6-isometries}Let $M$ be a minimal torus immersed
in $S^{3}$. If $nnt(M)\leq6$, then the group of isometries of $M$
is not trivial.
\end{thm}
\begin{proof}
Unless there is some nonzero $B\in so(4)$ for which $f_{B}=0$, in
which case Proposition (\ref{pro:f_B-vanish}) implies the existence
of a one-parameter group of isometries of $S^{3}$ which restrict
to isometries of $M$, then $nnt(M)\leq6$ implies that the span of
$\{u_{1},u_{2}\}$, $u_{1}:=a^{-\frac{3}{2}}W_{1}(a)=h_{0}$ and $u_{2}:=a^{-\frac{3}{2}}W_{2}(a)=h_{\frac{\pi}{2}}$,
will be contained in the span of $\{f_{B}|B\in so(4)\}$. Since then
$u_{1}=2f_{B}$ for some $B\in so(4)$, then $M$ defines a solution
$\phi$ of the system (\ref{eq:1}) associated with the matrix $B$
and with $\theta=0$. The condition $u_{2}=2f_{\tilde{B}}$ implies
by Remark(\ref{rem:h}) that $s=\tilde{\xi_{1}}$, for the system
(\ref{eq:2}) associated with the matrices $B,\,\tilde{B}$ and $\theta=0$.
As before, we will assume that $\xi_{1}(0,0)=s(0,0)=0$ and $r(0,0)=r_{0}$.
Define the function $f=s-\tilde{\xi}_{1}$. The hypothesis in the
theorem is equivalent to the condition that $f$ is identically zero,
in particular, $\tilde{\xi}_{1}(0,0)=0$, since $f(0,0)=0$. The theorem
is a consequence of the previous lemma and will follow by showing
that $\xi_{4}(0,0)=0$. A direct computation shows that\[
\frac{\partial f}{\partial u}=\e^{-r}\xi_{6}-\e^{r}\xi_{3}-\e^{-r}\tilde{\xi_{5}}-\e^{r}\tilde{\xi_{2}}\]
and\begin{eqnarray*}
\frac{\partial^{2}f}{\partial u^{2}} & = & \xi_{1}(-\e^{-r}\xi_{6}-\e^{r}\xi_{3}+\e^{-r}\tilde{\xi}_{5}-\e^{r}\tilde{\xi}_{2})\\
 &  & +\e^{-r}(-s\xi_{5}+\e^{r}\xi_{4})-\e^{r}(-s\xi_{2}-\e^{-r}\xi_{4})\\
 &  & -e^{-r}(-s\xi_{5}+\e^{r}\xi_{4})-\e^{r}(-s\xi_{2}-\e^{-r}\xi_{4})\\
 &  & -\e^{r}(s\tilde{\xi_{6}}-\e^{-r}\xi_{1})-\e^{r}(s\tilde{\xi}_{3}-\e^{r}\tilde{\xi}_{1})\\
 & = & \xi_{1}(-\e^{-r}\xi_{6}-\e^{r}\xi_{3}-\e^{-r}\tilde{\xi}_{5}-\e^{r}\tilde{\xi}_{2})\\
 &  & +s(-\e^{-r}\xi_{5}+\e^{r}\xi_{2}-\e^{-r}\tilde{\xi}_{6}-\e^{r}\tilde{\xi}_{3})\\
 &  & +2\xi_{4}+2\cosh(2r)\tilde{\xi_{1}}.\end{eqnarray*}
From the last equation, using the fact that $s(0,0)=\xi_{1}(0,0)=\tilde{\xi_{1}}(0,0)$
and $\frac{\partial^{2}f}{\partial u^{2}}=0$, we conclude that $\xi_{4}(0,0)=0$,
which implies, by the previous lemma, that $r(u,v)=r(-u,-v)$. To
finish the proof of the theorem, we notice that the function $A(u,v)=-(u,v)$
preserves the lattice in $\R^{2}$ given by the double-periodicity
of the function $\phi$ and therefore induces a function in the torus
$\tau(\R^{2})=M$, since the first fundamental form of $M$ in the
coordinates $u$ and $v$ is $c\e^{-2r}(du^{2}+dv^{2})$ where $c$
is a positive constant, then, this function from $M$ to $M$ induced
by $A$ is an isometry.
\end{proof}

\end{document}